\documentclass[12pt,reqno]{amsart}

\usepackage{amscd}
\usepackage{amssymb}
\usepackage{epsfig}
\usepackage{float}
\usepackage{url}
\usepackage[all]{xy}
\usepackage{color}

\newcommand{\N}{{\mathbb N}}
\newcommand{\Z}{{\mathbb Z}}
\newcommand{\Q}{{\mathbb Q}}
\newcommand{\C}{{\mathbb C}}
\newcommand{\R}{{\mathbb R}}
\renewcommand{\P}{{\mathbb P}}
\renewcommand{\H}{{\mathbb H}}

\newcommand{\OO}{{\mathcal O}}

\newcommand{\RR}{{\mathcal R}}

\newcommand{\www}{\widetilde}

\newcommand{\oooo}{\overline}
\newcommand{\uuuu}{\underline}

\newcommand{\mmm}{{\bf m}}

\newcommand{\ppp}{\partial}

\DeclareMathOperator{\Aut}{Aut}

\DeclareMathOperator{\id}{id}
\DeclareMathOperator{\Imm}{Im}

\DeclareMathOperator{\lcm}{lcm}
\DeclareMathOperator{\mult}{mult}
\DeclareMathOperator{\ord}{ord}
\DeclareMathOperator{\Ord}{Ord}

\DeclareMathOperator{\rank}{rank}

\DeclareMathOperator{\Rad}{Rad}

\DeclareMathOperator{\Stab}{Stab}

\begin{document}

\theoremstyle{plain}
\newtheorem{lemma}{Lemma}[section]
\newtheorem{definition/lemma}[lemma]{Definition/Lemma}
\newtheorem{theorem}[lemma]{Theorem}
\newtheorem{proposition}[lemma]{Proposition}
\newtheorem{corollary}[lemma]{Corollary}
\newtheorem{conjecture}[lemma]{Conjecture}
\newtheorem{conjectures}[lemma]{Conjectures}

\theoremstyle{definition}
\newtheorem{definition}[lemma]{Definition}
\newtheorem{withouttitle}[lemma]{}
\newtheorem{remark}[lemma]{Remark}
\newtheorem{remarks}[lemma]{Remarks}
\newtheorem{example}[lemma]{Example}
\newtheorem{examples}[lemma]{Examples}

\title
[$\mu$-constant monodromy groups for some singularities]
{$\mu$-constant monodromy groups and Torelli
results for marked singularities, for the
unimodal and some bimodal singularities}

\author{Falko Gau{\ss}\and Claus Hertling}

\address{Falko Gau{\ss}\\Universit\"at Mannheim\\ 
Lehrstuhl f\"ur Mathematik VI\\
Seminargeb\"aude A 5, 6\\
68131 Mannheim, Germany}

\email{gauss@math.uni-mannheim.de}

\address{Claus Hertling\\Universit\"at Mannheim\\ 
Lehrstuhl f\"ur Mathematik VI\\
Seminargeb\"aude A 5, 6\\
68131 Mannheim, Germany}

\email{hertling@math.uni-mannheim.de}

\thanks{This work was supported by the DFG grant He2287/4-1
(SISYPH)}



\date{April 11, 2016}

\begin{abstract}
This paper is a sequel to \cite{He7}. There a notion of
marking of isolated hypersurface singularities was defined, 
and a moduli space $M_\mu^{mar}$ 
for marked singularities in one $\mu$-homotopy class 
of isolated hypersurface singularities was established.
One can consider it as a global $\mu$-constant stratum
or as a Teichm\"uller space for singularities.
It comes together with a $\mu$-constant monodromy group 
$G^{mar}\subset G_\Z$. Here $G_\Z$ is the group of automorphisms
of a Milnor lattice which respect the Seifert form.
It was conjectured that $M_\mu^{mar}$ is connected.
This is equivalent to $G^{mar}= G_\Z$.
Also Torelli type conjectures were formulated.
All conjectures were proved for the simple singularities and
22 of the exceptional unimodal and bimodal singularities.
In this paper the conjectures are proved for the remaining
unimodal singularities and the remaining exceptional bimodal
singularities.
\end{abstract}

\maketitle

\tableofcontents

\setcounter{section}{0}

\section{Introduction}\label{c1}

\noindent
This paper is a sequel to \cite{He7}.
That paper studied local objects, namely 
holomorphic functions germs
$f:(\C^{n+1},0)\to (\C,0)$ with an isolated singularity
at 0 (short: singularity), from a global perspective.

There a notion of marking of a singularity was defined.
One has to fix one singularity $f_0$, which serves
as reference singularity. Then a marked singularity is a pair
$(f,\pm\rho)$ where $f$ is in the $\mu$-homotopy class of $f_0$
and $\rho:(Ml(f),L)\to (Ml(f_0),L)$ is an isomorphism.
Here $Ml(f)$ is the Milnor lattice of $f$, and $L$ is the 
Seifert form on it (the definitions are recalled in section 
\ref{c2}). The group $G_\Z(f_0):=\textup{Aut}(Ml(f_0),L)$
will be important, too.

A moduli space $M_\mu^{mar}$ of right equivalence classes 
of marked singularities in one $\mu$-homotopy class was established.
One can consider it as a global $\mu$-constant stratum
or as a Teichm\"uller space for singularities.
The group $G_\Z(f_0)$ acts properly discontinuously on it, 
and the quotient is the
moduli space $M_\mu$ of right equivalence classes of unmarked
singularities from \cite[chapter 13]{He6}.
The {\it $\mu$-constant monodromy group}
$G^{mar}(f_0)\subset G_\Z(f_0)$ is the subgroup of 
automorphisms which map the topological component
$(M^{mar}_\mu(f_0))^0$ which contains $[(f_0,\pm\id)]$
to itself. It can also be constructed as the group of all
automorphisms which can be realized modulo $\pm\id$ as
transversal monodromies of $\mu$-constant families.

\begin{conjecture}\label{t1.1} \cite[conjecture 3.2 (a)]{He7}
$M^{mar}_\mu$ is connected. Equivalently: $G^{mar}=G_\Z$.
\end{conjecture}

Roughly this conjecture says that all abstract
automorphisms come from geometry, from coordinate changes
and $\mu$-constant families. 

Also Torelli type conjectures are formulated in \cite{He7}.
Any singularity $f$ comes equipped with its {\it Brieskorn lattice}
$H_0''(f)$, and any marked singularity $(f,\pm\rho)$ 
comes equipped with a {\it marked Brieskorn lattice} $BL(f,\pm\rho)$.
The Gau{\ss}-Manin connection for singularities and the
Brieskorn lattices had been introduced 1970 by Brieskorn
and had been studied since then. The second author has
a long-going project on Torelli type conjectures around them.

In \cite{He4} a classifying space $D_{BL}(f_0)$ for
marked Brieskorn lattices was constructed. It is especially
a complex manifold, and $G_\Z(f_0)$ acts properly discontinuously
on it. The quotient $D_{BL}/G_{\Z}$ is a space of isomorphism
classes of Brieskorn lattices. There is a natural holomorphic 
period map 
$$BL:M^{mar}_\mu(f_0)\to D_{BL}(f_0).$$
It is an immersion \cite[theorem 12.8]{He6}
(this refines a weaker result in \cite{SaM}).
And it is $G_\Z$-equivariant. 

\begin{conjecture}\label{t1.2}
(a) \cite[conjecture 5.3]{He7}
$BL$ is an embedding.

(b) \cite[conjecture 5.4]{He7}, \cite[conjecture 12.7]{He6},
\cite[Kap. 2 d)]{He1}) $LBL: M_\mu=M^{mar}_\mu/G_\Z\to D_{BL}/G_\Z$
is an embedding.
\end{conjecture}

Part (b) says that the right equivalence class of $f$
is determined by the isomorphism class of $H_0''(f)$.
Part (a) would imply part (b). Both are global Torelli type
conjectures. Part (b) was proved in \cite{He1} for all
simple and unimodal singularities and almost all bimodal
singularities, all except three subseries of the eight
bimodal series. Therefore for the proof of part (a)
in these cases it remains mainly to control $G_\Z$ well.
But that is surprisingly difficult.

In \cite{He7} the conjectures \ref{t1.1} and \ref{t1.2}
were proved for all simple and those 22 of the 28 exceptional
unimodal and bimodal singularities, where all eigenvalues of
the monodromy have multiplicity one.
In this paper they will be proved for the remaining unimodal
and exceptional bimodal singularities, that means,
for the simple elliptic and the hyperbolic singularities
and for those 6 of the 28 exceptional unimodal and bimodal
singularities which had not been treated in \cite{He7}.

A priori, logically conjecture \ref{t1.1} comes before
conjecture \ref{t1.2}. But the results in \cite{He1} are
more concrete and give already some information about 
the action of $G_\Z$ on $D_{BL}$ and $M^{mar}_\mu$.
Anyway, the main remaining work is a good control of the 
groups $G_\Z$.
That presents some unexpected difficulties. For example,
we need two surprising generalizations of 
the number theoretic fact 
$\Z[e^{2\pi i/a}]\cap S^1=\{\pm e^{2\pi ik/a}\, |\, k\in\Z\}$
for $a\in\N$, one is lemma \ref{t2.5}, the other
is related to $U_{16}$, see remark \ref{t4.3}.

The groups $G_\Z(f_0)$ will be calculated 
in the sections \ref{c3} and \ref{c4}.
Section \ref{c2} collects well known background material 
on the topology of singularities. But it contains also an
algebraic lemma \ref{t2.5} about automorphisms of 
monodromy modules. Section \ref{c5} collects notions and
results from \cite{He7} on marked singularities, the 
moduli spaces $M^{mar}_\mu(f_0)$, the groups 
$G^{mar}(f_0)$ and the Torelli type conjectures.
The sections \ref{c6} and \ref{c7} give the proofs
of the conjectures \ref{t1.1} and \ref{t1.2} in the cases
considered.
Section \ref{c8} is motivated by the paper \cite{MS} of Milanov
and Shen and complements their results on (transversal) 
monodromy groups for certain families of simple elliptic singularities.
The three principal congruence subgroups $\Gamma(3),\Gamma(4)$
and $\Gamma(6)$ which turn up in \cite{MS} for certain families
are shown to turn up also in the biggest possible families.

\section{Review on the topology of isolated
hypersurface singularities}\label{c2}

\noindent
First, we recall some classical facts and fix some notations.
An {\it isolated hypersurface singularity} (short: {\it singularity})
is a holomorphic function germ $f:(\C^{n+1},0)\to (\C,0)$ with an isolated 
singularity at $0$. Its {\it Milnor number} 
$$\mu:=\dim\OO_{\C^{n+1},0}/(\frac{\ppp f}{\ppp x_i})$$ 
is finite. 
For the following notions and facts compare \cite{AGV2}, 
\cite{Eb2} and (for the notion of an unfolding) \cite{AGV1}.
A {\it good representative} of $f$ has to be defined with some 
care \cite{Mi}\cite{AGV2}\cite{Eb2}. It is $f:Y\to T$
with $Y\subset\C^{n+1}$ a suitable small neighborhood of 0 and 
$T\subset\C$ a small disk around 0.
Then $f:Y'\to T'$ with $Y'=Y-f^{-1}(0)$ and 
$T'=T-\{0\}$ is a locally trivial $C^\infty$-fibration,
the  {\it Milnor fibration}. Each fiber has the
homotopy type of a bouquet of $\mu$ $n$-spheres \cite{Mi}.

Therefore the (reduced for $n=0$) middle homology groups are {}\\{}
$H_n^{(red)}(f^{-1}(\tau),\Z) \cong \Z^\mu$ for $\tau\in T'$.
Each comes equipped with an intersection form $I$, 
which is a datum of one fiber,
a monodromy $M_h$ and a Seifert form $L$, which come from the 
Milnor fibration,
see \cite[I.2.3]{AGV2} for their definitions 
(for the Seifert form, there are several
conventions in the literature, we follow \cite{AGV2}). 
$M_h$ is a quasiunipotent automorphism, $I$ and $L$ are 
bilinear forms with values in $\Z$,
$I$ is $(-1)^n$-symmetric, and $L$ is unimodular. $
L$ determines $M_h$ and $I$ because of the formulas
\cite[I.2.3]{AGV2}
\begin{eqnarray}\label{2.1}
L(M_ha,b)&=&(-1)^{n+1}L(b,a),\\ \label{2.2}
I(a,b)&=&-L(a,b)+(-1)^{n+1}L(b,a).
\end{eqnarray}
The Milnor lattices $H_n(f^{-1}(\tau),\Z)$ for all Milnor fibrations
$f:Y'\to T'$ and then all $\tau\in\R_{>0}\cap T'$ are canonically isomorphic,
and the isomorphisms respect $M_h$, $I$ and $L$. 
This follows from Lemma 2.2 in \cite{LR}. 
These lattices are identified and called {\it Milnor lattice} $Ml(f)$.
The group $G_\Z$ is 
\begin{eqnarray}\label{2.3}
G_\Z=G_\Z(f):= \Aut(Ml(f),L)=\Aut(Ml(f),M_h,I,L),
\end{eqnarray}
the second equality is true because $L$ determines $M_h$ and $I$.
We will use the notation $Ml(f)_\C:=Ml(f)\otimes_\Z \C$,
and analogously for other rings $R$ with $\Z\subset R\subset \C$,
and the notations
\begin{eqnarray*}
Ml(f)_\lambda&:=&\ker((M_h-\lambda\id)^\mu:Ml(f)_\C\to Ml(f)_\C)
\subset Ml(f)_\C,\\
Ml(f)_{1,\Z}&:=& Ml(f)_1\cap Ml(f)\subset Ml(f),\\
Ml(f)_{\neq 1}&:=&\bigoplus_{\lambda\neq 1}Ml(f)_\lambda \subset 
Ml(f)_\C,\\
Ml(f)_{\neq 1,\Z}&:=& Ml(f)_{\neq 1}\cap Ml(f)\subset Ml(f).
\end{eqnarray*}

The formulas \eqref{2.1} and \eqref{2.2} show
$I(a,b)= L((M_h-\id)a,b)$. Therefore the eigenspace with eigenvalue
1 of $M_h$ is the radical $\Rad(I)\subset Ml(f)$ of $I$. 
By \eqref{2.2} $L$ is 
$(-1)^{n+1}$-symmetric on the radical of $I$.

In the case of a curve singularity ($n=1$) with $r$ branches, 
$f=\prod_{j=1}^r f_j$,
the radical of $I$ is a $\Z$-lattice of rank $r-1$, 
and it is generated by the classes $l_j\in Ml(f)$ which are obtained
by pushing the (correctly oriented) 
cycles $\partial Y\cap f_j^{-1}(0)$ from
the boundary of the fiber $f^{-1}(0)$ to the boundary of the
fiber $f^{-1}(\tau)$. Then
\begin{eqnarray}\label{2.4}
l_1+...+l_r&=&0,\\
L(l_i,l_j)&=&\textup{intersection number of }(f_i,f_j)\quad
\textup{for }i\neq j,\nonumber\\
\textup{so }L(l_i,l_j)&>&0\quad\textup{ for }i\neq j,\label{2.5}\\
l_1,...,\widehat{l_j},...,l_r&&\textup{is a }\Z\textup{-basis
of }\Rad(I)\textup{ for any }j\label{2.6}.
\end{eqnarray}
Kaenders proved the following result using Selling reduction.
It will be useful for the calculation of $\Aut(\Rad(I),L)$,
because it implies that any automorphism of $(\Rad(I),L)$ 
maps the set $\{l_1,...,l_r\}$ to itself or to minus itself.

\begin{theorem}\label{t2.1}\cite{Ka}
In the case of a curve singularity as above, the set 
$\{l_1,...,l_r\}$ is up to a common sign uniquely determined
by the properties \eqref{2.4}, \eqref{2.5} and \eqref{2.6}.
So it is up to a common sign determined by the pair $(\Rad(I),L)$.
Furthermore, $L$ is negative definite on $\Rad(I)$.
\end{theorem}

\begin{examples}\label{t2.2}
In the following three examples $\uuuu{l}=(l_1,...,l_r)$,
and $L(\uuuu{l}^t,\uuuu{l})=(L(l_i,l_j))$.
Then $(L(l_i,l_j))_{1\leq i,j\leq r-1}$ is the matrix
of $L$ on $\Rad(I)$ with respect to the $\Z$-basis $l_1,...,l_{r-1}$.
The examples (ii) and (iii) will be useful in section \ref{c4},
the example (i) is an alternative to a calculation in the
proof of theorem 8.4 in \cite{He7}.

\medskip
(i) $D_{2m}:\ f=x^{2m-1}+xy^2=x(x^{m-1}-iy)(x^{m-1}+iy),$
\begin{eqnarray*}
L(\uuuu{l}^t,\uuuu{l})&=& 
\begin{pmatrix}-2&1&1\\ 1&-m&m-1\\ 1&m-1&-m\end{pmatrix},\quad
(L(l_i,l_j))_{1\leq i,j\leq 2}=
\begin{pmatrix}-2&1\\ 1&-m\end{pmatrix}.
\end{eqnarray*}
Obviously $|\Aut(\Rad(I),L)|=12$ in the case $m=2$, and
$|\Aut(\Rad(I),L)|=4$ in the cases $m\geq 3$.

\medskip
(ii) $Z_{12}:\ f=x^3y+xy^4=xy(x^2+y^3),$
\begin{eqnarray*}
L(\uuuu{l}^t,\uuuu{l})&=& 
\begin{pmatrix}-4&1&3\\ 1&-3&2\\ 3&2&-5\end{pmatrix},\qquad
(L(l_i,l_j))_{1\leq i,j\leq 2}=
\begin{pmatrix}-4&1\\ 1&-3\end{pmatrix}.
\end{eqnarray*}
Obviously $\Aut(\Rad(I),L)=\{\pm\id\}$. 

\medskip
(iii) $Z_{18}:\ f=x^3y+xy^6=xy(x^2+y^5),$
\begin{eqnarray*}
L(\uuuu{l}^t,\uuuu{l})&=& 
\begin{pmatrix}-6&1&5\\ 1&-3&2\\ 5&2&-7\end{pmatrix},\qquad
(L(l_i,l_j))_{1\leq i,j\leq 2}=
\begin{pmatrix}-6&1\\ 1&-3\end{pmatrix}.
\end{eqnarray*}
Obviously $\Aut(\Rad(I),L)=\{\pm\id\}$. 
\end{examples}

Finally, in section \ref{c3}, distinguished bases will be used.
Again, good references for them are \cite{AGV2} and \cite{Eb2}.
We sketch their construction and properties.

One can choose a {\it universal unfolding} of $f$,
a {\it good representative} $F$ of it with base space 
$M\subset\C^\mu$,
and a generic parameter $t\in M$. Then $F_t:Y_t\to T$
with $T\subset \C$ the same disk as above and $Y_t\subset\C^{n+1}$
is a {\it morsification} of $f$. It has $\mu$ $A_1$-singularities,
and their critical values $u_1,...,u_\mu\in T$ 
are pairwise different. Their numbering is also a choice. 
Now choose a value $\tau\in T\cap\R_{>0}-\{u_1,...,u_\mu\}$ and
a {\it distinguished system of paths}. That is
a system of $\mu$ paths $\gamma_j$, $j=1,...,\mu$, from
$u_j$ to $\tau$ which do not intersect except at $\tau$
and which arrive at $\tau$ in clockwise order.
Finally, shift from the $A_1$ singularity above each value $u_j$
the (up to sign unique) vanishing cycle along $\gamma_j$
to the Milnor fiber $Ml(f)=H_n(f^{-1}(\tau),\Z)$,
and call the image $\delta_j$. 

The tuple $(\delta_1,...,\delta_\mu)$ forms a $\Z$-basis of 
$Ml(f)$. All such bases are called {\it distinguished bases}.
They form one orbit of an action of a semidirect product 
$\textup{Br}_\mu\ltimes \{\pm 1\}^\mu$. Here $\textup{Br}_\mu$ is the braid
group with $\mu$ strings, see \cite{AGV2} or \cite{Eb2} for its action. The {\it sign change group} $\{\pm 1\}^\mu$ acts simply
by changing the signs of the entries of the tuples
$(\delta_1,...,\delta_\mu)$. 
The members of the distinguished bases are called 
{\it vanishing cycles}.

The {\it Stokes matrix} $S$ of a distinguished basis 
is defined as the upper triangular matrix in $M(\mu\times\mu,\Z)$
with 1's on the diagonal and with
$$S_{ij}:=(-1)^{n(n+1)/2}\cdot I(\delta_i,\delta_j) \quad
\textup{for all }i,j\textup{ with }i<j.$$
The {\it Coxeter-Dynkin diagram} of a distinguished basis encodes
$S$ in a geometric way. It has $\mu$ vertices which are numbered
from 1 to $\mu$. Between two vertices $i$ and $j$ with $i<j$
one draws

\begin{tabular}{ll}
no edge & if $S_{ij}=0$, \\
$|S_{ij}|$ edges & if $S_{ij}<0$, \\
$S_{ij}$ dotted edges & if $S_{ij}>0$. 
\end{tabular}

Coxeter-Dynkin diagrams of many singularities were calculcated
by A'Campo, Ebeling, Gabrielov and Gusein-Zade. Some of them
can be found in \cite{Ga}, \cite{Eb1} and \cite{Eb2}.

\begin{example}\label{t2.3}
The hyperbolic singularities of type $T_{pqr}$ with 
$\frac{1}{p}+\frac{1}{q}+\frac{1}{r}<1$ and the
simple elliptic singularities of types $\www E_6=T_{333},
\www E_7=T_{442}$ and $\www E_8=T_{632}$ have 
distinguished bases with the Coxeter-Dynkin diagrams
in figure 1 \cite{Ga}.

\begin{figure}[h]
\begin{center}
\includegraphics[width=0.8\textwidth]{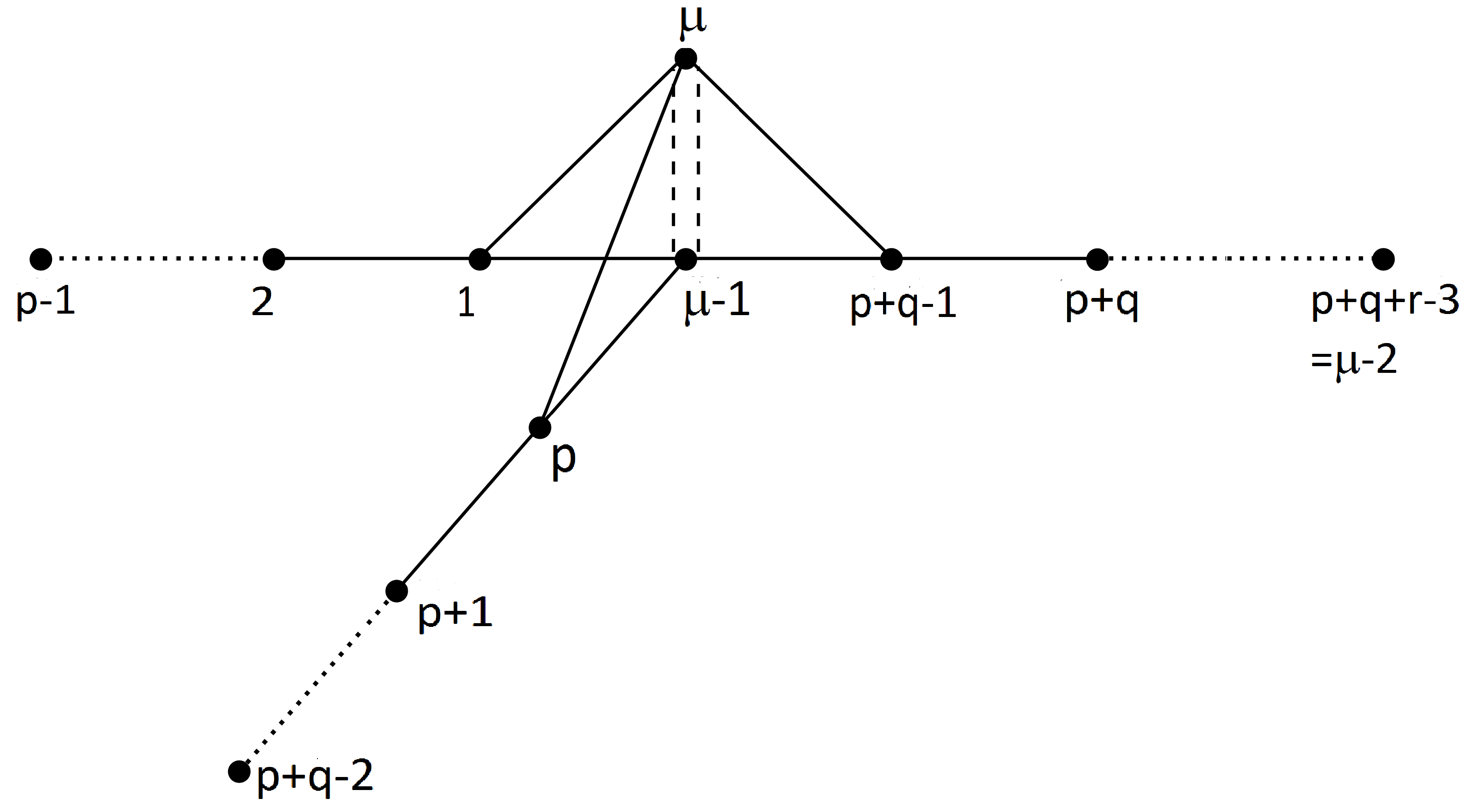} 
\end{center}
\caption{A Coxeter-Dynkin diagram of the singularities of type 
$T_{pqr}$}\label{pic-1}
\end{figure}

\end{example}

The Picard-Lefschetz transformation on $Ml(f)$ 
of a vanishing cycle $\delta$ is
\begin{eqnarray}\label{2.7}
s_{\delta}(b)&:=&b-(-1)^{n(n+1)/2}\cdot I(\delta,b)\cdot \delta.
\end{eqnarray}
The monodromy $M_h$ is
\begin{eqnarray}\label{2.8}
M_h &=& s_{\delta_1}\circ ...\circ s_{\delta_\mu}
\end{eqnarray}
for any distinguished basis $\uuuu{\delta}=(\delta_1,...,\delta_\mu)$.
The matrices of the monodromy, Seifert form and intersection form
with respect to a distinguished basis 
$\uuuu{\delta}$ are given by the following formulas.
\begin{eqnarray}\label{2.9}
M_h(\uuuu{\delta})&=&\uuuu{\delta}\cdot (-1)^{n+1}\cdot S^{-1}S^t,\\
I(\uuuu{\delta}^t,\uuuu{\delta})&=&(-1)^{n(n+1)/2}\cdot 
(S+(-1)^nS^t),\label{2.10}\\
L(\uuuu{\delta}^t,\uuuu{\delta})&=&(-1)^{(n+1)(n+2)/2}\cdot S^t.\label{2.11}
\end{eqnarray}

\begin{remark}\label{t2.4}
The Stokes matrix $S$ of a distinguished basis is 
related to the matrix $V$ in \cite[Korollar 5.3 (i)]{Eb2}
by the formula $V=L(\uuuu{\delta}^t,\uuuu{\delta})=
(-1)^{(n+1)(n+2)/2}\cdot S^t$. Thus $V$ is
{\it lower} triangular, not {\it upper} triangular,
contrary to the claim in \cite[Korollar 5.3 (i)]{Eb2}.
The matrix of $\textup{Var}^{-1}$ for a distinguished 
basis $\uuuu\delta$ and its dual basis $\uuuu\delta^*$
is $(-1)^{(n+1)(n+2)/2}S$, namely 
$$\textup{Var}^{-1}(\uuuu\delta)
=\uuuu\delta^*\cdot (-1)^{(n+1)(n+2)/2}\cdot S.$$
But then the matrix for the Seifert form $L$ with
$L(a,b)=(\textup{Var}^{-1}(a))(b)$ is
$$L(\uuuu{\delta}^t,\uuuu{\delta})
=\textup{Var}^{-1}(\uuuu{\delta}^t)(\uuuu{\delta})
=(-1)^{(n+1)(n+2)/2}S^t.$$
In \cite[I.2.5]{AGV2} this problem does not arise
because there the matrix of a bilinear form with respect to
a basis is the transpose of the usual matrix \cite[page 45]{AGV2}.
This is applied to the matrices of $I$ and $L$.
\end{remark}

A result of Thom and Sebastiani
compares the Milnor lattices and monodromies of 
the singularities $f=f(x_0,...,x_n),g=g(y_0,...,y_m)$ and
$f+g=f(x_0,...,x_n)+g(x_{n+1},...,x_{n+m+1})$.
There are extensions by Deligne for the Seifert form and 
by Gabrielov for distinguished bases. 
All results can be found in \cite[I.2.7]{AGV2}.
They are restated here.
There is a canonical isomorphism
\begin{eqnarray}\label{2.12}
\Phi:Ml(f+g)&\stackrel{\cong}{\longrightarrow} &Ml(f)\otimes Ml(g),\\
\textup{with } M_h(f+g)&\cong & M_h(f)\otimes M_h(g) \label{2.13}\\
\textup{and } 
L(f+g)&\cong& (-1)^{(n+1)(m+1)}\cdot L(f)\otimes L(g).\label{2.14}
\end{eqnarray}
If $\uuuu{\delta}=(\delta_1,...,\delta_{\mu(f)})$
and $\uuuu{\gamma}=(\gamma_1,...,\gamma_{\mu(g)})$ are
distinguished bases of $f$ and $g$  with Stokes matrices
$S(f)$ and $S(g)$, then 
$$\Phi^{-1}(\delta_1\otimes \gamma_1,...,
\delta_1\otimes \gamma_{\mu(g)},
\delta_2\otimes \gamma_1,...,
\delta_2\otimes \gamma_{\mu(g)},
...,
\delta_{\mu(f)}\otimes \gamma_1,...,
\delta_{\mu(f)}\otimes \gamma_{\mu(g)})$$
is a distinguished basis of $Ml(f+g)$,
that means, one takes the vanishing cycles 
$\Phi^{-1}(\delta_i\otimes \gamma_j)$ in the lexicographic order.
Then by \eqref{2.11} and \eqref{2.14}, the matrix 
\begin{eqnarray}\label{2.15}
S(f+g)=S(f)\otimes S(g)
\end{eqnarray}
(where the tensor product is defined
so that it fits to the lexicographic order) 
is the Stokes matrix of this distinguished basis.

In the special case $g=x_{n+1}^2$,
the function germ $f+g=f(x_0,...,x_n)+x_{n+1}^2\in \OO_{\C^{n+2},0}$
is called {\it stabilization} or {\it suspension} of $f$. 
As there are only two isomorphisms $Ml(x_{n+1}^2)\to\Z$, 
and they differ by a sign, there are two equally canonical
isomorphisms $Ml(f)\to Ml(f+x_{n+1}^2)$,
and they differ just by a sign. 
Therefore automorphisms and bilinear forms on $Ml(f)$ 
can be identified with automorphisms and bilinear forms on 
$Ml(f+x_{n+1}^2)$. In this sense
\begin{eqnarray}\label{2.16}
L(f+x_{n+1}^2) = (-1)^n\cdot L(f)\quad\textup{and}\quad 
M_h(f+x_{n+1}^2)= - M_h(f)
\end{eqnarray}
\cite[I.2.7]{AGV2}, and $G_\Z(f+x_{n+1}^2)=G_\Z(f)$.
The Stokes matrix $S$ does not change under stabilization.

The following algebraic lemma from \cite{He7} will be
very useful in the sections \ref{c3} and \ref{c4}. 
It can be seen as a generalization of the number theoretic fact
$\Z[e^{2\pi i/a}]\cap S^1=\{\pm e^{2\pi ik/a}\, |\, k\in\Z\}$.

\begin{lemma}\label{t2.5} \cite[lemma 8.2]{He7}
Let $H$ be a free $\Z$-module of finite rank $\mu$, and $H_\C:=H\otimes_\Z\C$.
Let $M_h:H\to H$ be an automorphism of finite order, called monodromy,
with three properties:
\begin{list}{}{}
\item[(i)] 
Each eigenvalue has multiplicity $1$. \\
Denote $H_\lambda:=\ker(M_h-\lambda\cdot\id:H_\C\to H_\C).$
\item[(ii)]
Denote 
$\Ord:=\{\ord\lambda\, |\, \lambda\textup{ eigenvalue of }M_h\}
\subset \Z_{\geq_1}$. 
There exist four sequences $(m_i)_{i=1,...,|\Ord|}$, $(j(i))_{i=2,...,|\Ord|}$, 
$(p_i)_{i=2,...,|\Ord|}$, $(k_i)_{i=2,...,|\Ord|}$ of numbers in $\Z_{\geq 1}$
and two numbers $i_1,i_2\in\Z_{\geq 1}$ with $i_1\leq i_2\leq |\Ord|$ and 
with the properties: \\
$\Ord=\{m_1,...,m_{|\Ord|}\}$,\\ 
$p_i$ is a prime number, $p_i=2$ for $i_1+1\leq i\leq i_2$, $p_i\geq 3$ else, \\
$j(i)=i-1$ for $i_1+1\leq i\leq i_2$, $j(i)<i$ else,
$$m_i=m_{j(i)}/p_i^{k_i}.$$
\item[(iii)]
A cyclic generator $a_1\in H$ exists, that means,
$$H=\bigoplus_{i=0}^{\mu-1}\Z\cdot M_h^i(a_1).$$
\end{list}
Finally, let $I$ be an $M_h$-invariant nondegenerate bilinear form
(not necessarily $(\pm 1)$-symmetric) on $\bigoplus_{\lambda\neq \pm 1}H_\lambda$
with values in $\C$. Then
$$\Aut(H,M_h,I)=\{\pm M_h^k\, |\, k\in \Z\}.$$
\end{lemma}

\section{The group $G_\Z$ for the simple elliptic
and the hyperbolic singularities}\label{c3}

\noindent
The simple elliptic and the hyperbolic singularities are
1-parameter families of singularities, which had been
classified by Arnold \cite{AGV1}. For each
triple $(p,q,r)\in\N_{\geq 2}^3$ with $p\geq q\geq r$ and 
$\kappa:=\frac{1}{p}+\frac{1}{q}+\frac{1}{r}\leq 1$
one has one family, denoted $T_{pqr}$.
The hyperbolic singularities are those with $\kappa<1$,
the simple elliptic are those with $\kappa=1$.
For the three families of simple elliptic singularities
also other symbols are used,
$T_{333}=\www E_6,\ T_{442}=\www E_7,\ T_{632}=\www E_8$.
The singularities of types $T_{pqr}$ with $r=2$ 
exist as curve singularities, all others
as surface singularities. Normal forms will be discussed
in section \ref{c6}. Here for each family the group 
$G_\Z=\Aut(Ml(f),L)$ will be analyzed. The result in 
Theorem \ref{t3.1} is completely explicit in the case $\kappa<1$
and partly explicit in the case $\kappa=1$.
The whole section is devoted to its proof.
Besides 
$$\kappa:=\frac{1}{p}+\frac{1}{q}+\frac{1}{r},\quad
\textup{also}\quad \chi:=\textup{lcm}(p,q,r)\in\N$$ will be used.
A singularity of type $T_{pqr}$ has Milnor number $\mu=p+q+r-1$,
and its monodromy has the characteristic polynomial
$$\frac{t^p-1}{t-1}\cdot\frac{t^q-1}{t-1}\cdot\frac{t^r-1}{t-1}
\cdot (t-1)^2.$$

\begin{theorem}\label{t3.1}
Consider a surface singularity $f$ of type $T_{pqr}$ 
(with $\kappa\leq 1$) with Milnor lattice $Ml(f)$,
monodromy $M_h$, intersection form $I$ and Seifert form $L$.

\medskip
(a) Then $\dim Ml(f)_1=2,\ \rank\Rad(I)=1$ if $\kappa<1$
and $=2$ if $\kappa=1$. Choose a $\Z$-basis $b_1,b_2$ of
$Ml(f)_{1,\Z}$ with $b_1\in \Rad(I)$ and 
$L(b_1,b_2)\leq 0$. Then 
\begin{eqnarray}\label{3.1}
L(b_1,b_2)=-\chi\quad\textup{and}\quad 
M_hb_2=b_2+\chi(\kappa-1)\cdot b_1.
\end{eqnarray}

\medskip
(b) The restriction map $G_\Z\to \Aut(Ml(f)_{1,\Z},L)$ is surjective.
Denote by $T\in \Aut(Ml(f)_{1,\Z})$ the automorphism
with $T(b_1)=b_1$ and $T(b_2)=b_2+b_1$. Denote $\uuuu{b}:=(b_1,b_2)$.
\begin{eqnarray}\label{3.2}
\Aut(Ml(f)_{1,\Z},L)&=&\{\uuuu{b}\mapsto \uuuu{b}\cdot A\, |\, 
A\in SL(2,\Z)\} \\
&\cong& SL(2,\Z) \quad\textup{if}\quad \kappa=1,\nonumber\\
\Aut(Ml(f)_{1,\Z},L)&=& \{\pm T^k\, |\, k\in\Z\}
\quad\textup{if}\quad\kappa<1.\label{3.3}
\end{eqnarray}

\medskip
(c) The group $G_\Z$ for $\kappa<1$ and the subgroup
$\{g\in G_\Z\, |\, g(b_1)=\pm b_1\}\subset G_\Z$ for $\kappa=1$
will be described explicitly, except for the part
$U_2$, see below.
There is a monodromy invariant decomposition 
\begin{eqnarray}\label{3.4}
Ml(f)_{\neq 1}=Ml^{(1)}_\C\oplus Ml^{(2)}_\C\oplus Ml^{(3)}_\C
\end{eqnarray}
such that the characteristic polynomial of $M_h|_{Ml^{(j)}}$
is 
\begin{eqnarray}\label{3.5}
\frac{t^p-1}{t-1},\ \frac{t^q-1}{t-1},\ \frac{t^r-1}{t-1}
\quad\textup{ for }j=1,2,3
\end{eqnarray}
and such that the following holds.
\begin{eqnarray}\label{3.6}
\left. \begin{array}{r} 
G_\Z\quad\textup{for }\kappa<1\\
\{g\in G_\Z\, |\, g(b_1)=\pm b_1\}\quad\textup{for }\kappa=1
\end{array}\right\} =(U_1\rtimes U_2)\times \{\pm\id\},
\end{eqnarray}
where $U_1$ is the infinite subgroup of $G_\Z$
\begin{eqnarray}\label{3.7}
U_1 &=&\{T^\delta\times (M_h|_{Ml^{(1)}_\C})^\alpha
\times (M_h|_{Ml^{(2)}_\C})^\beta
\times (M_h|_{Ml^{(3)}_\C})^\gamma \, |\\
&& (\delta,\alpha,\beta,\gamma)\in\Z\times \Z_p\times \Z_q
\times \Z_r\textup{ with }
\frac{\alpha}{p}+\frac{\beta}{q}+\frac{\gamma}{r}
\equiv \frac{\delta}{\chi}\mod 1\}\nonumber
\end{eqnarray}
and where $U_2$ is a finite subgroup of $G_\Z$ with 
\begin{eqnarray}\label{3.8}
U_2\left\{\begin{array}{ll} = \{\id\}\quad&\textup{if }p>q>r,\\
\cong S_2\quad&\textup{if }p=q>r\textup{ or }p>q=r,\\
\cong S_3\quad&\textup{if }p=q=r.\end{array}\right.
\end{eqnarray}
which consists of certain 
automorphisms which act trivially on $Ml(f)_1$
and which permute those of the subspaces $Ml^{(j)}_\C$
which have equal dimension.
\end{theorem}

{\bf Proof:}
Choose a distinguished basis 
$\uuuu{\delta}=(\delta_1,...,\delta_\mu)$ with the 
Coxeter-Dynkin diagram in example \ref{t2.3}.
Then the monodromy matrix $M_M$ with 
$M_h(\uuuu{\delta})=\uuuu{\delta}\cdot M_M$ can be calculated
either with \eqref{2.8} or with \eqref{2.9}.
It had been calculated in \cite{He1} with \eqref{2.8}, and it is
(here all not specified entries are 0)
\begin{eqnarray}\label{3.9}
M_M&=& 
\begin{pmatrix}
M_1 & & & M_8\\
 & M_2 & & M_9\\
 & & M_3 & M_{10}\\
M_5 & M_6 & M_7 & M_4\end{pmatrix}
\end{eqnarray}
with the following blocks,
\begin{eqnarray*}
M_1 &=&\begin{pmatrix}
0 & & & -1 \\
1 & \ddots & & -1 \\
  & \ddots & 0 & -1 \\
  & & 1 & -1 \end{pmatrix}\in M((p-1)\times (p-1),\Z),\\
M_2 &\in& M((q-1)\times (q-1),\Z)\textup{ and }\\
M_3 &\in& M((r-1)\times (r-1),\Z)\textup{ are defined analogously},\\
M_4&=& \begin{pmatrix}3&2\\ -2&-1\end{pmatrix},
\end{eqnarray*}
and $M_5,M_6,M_7,M_8,M_9,M_{10}$ are of suitable sizes with
all entries except the following being $0$,
\begin{eqnarray*}
&&(M_5)_{11}=(M_6)_{11}=(M_7)_{11}=-1,\ \ 
(M_5)_{21}=(M_6)_{21}=(M_7)_{21}=1, \\
&&(M_8)_{11}=(M_8)_{12}=(M_9)_{11}=(M_9)_{12}
=(M_{10})_{11}=(M_{10})_{12}=1.
\end{eqnarray*}

Define
\begin{eqnarray}\label{3.10}
\www b_1&:=& \delta_{\mu-1}-\delta_\mu,\\
\www b_2&:=& \chi\cdot\left(\sum_{i=1}^{p-1}\frac{p-i}{p}\delta_i
+\sum_{i=1}^{q-1}\frac{q-i}{q}\delta_{p-1+i}\right. \label{3.11}\\
&&\left. \hspace*{2cm}
+\sum_{i=1}^{r-1}\frac{r-i}{r}\delta_{p+q-2+i}+\delta_{\mu-1}\right).
\nonumber
\end{eqnarray}
Then one calculates 
\begin{eqnarray}\label{3.12}
M_h(\www b_1)=\www b_1,\quad
M_h(\www b_2)=\www b_2+\chi(\kappa-1)\cdot \www b_1,
\end{eqnarray}
and with \eqref{2.11}, which is here
$L(\uuuu{\delta}^t,\uuuu{\delta})=S^t$,
\begin{eqnarray}
\begin{pmatrix}
L(\www b_1,\www b_1)& L(\www b_1,\www b_2)\\
L(\www b_2,\www b_1)& L(\www b_2,\www b_2)\end{pmatrix}
=\begin{pmatrix}0&-\chi\\\chi &\frac{\chi^2}{2}(\kappa-1)
\end{pmatrix}.\label{3.13}
\end{eqnarray}
By \eqref{3.12}, $\www b_1,\www b_2$ is a $\Q$-basis of
$Ml(f)_{1,\Q}$ and $M_h$ is on $Ml(f)_1$ semisimple if $\kappa=1$
and it has a $2\times 2$ Jordan block if $\kappa<1$ 
(of course, this is well known).
From the coefficients one sees that $\www b_1,\www b_2$ is 
a $\Z$-basis of $Ml(f)_{1,\Z}$. Here it is important
that the coefficients of $\www b_2$ have greatest common divisor 1.
As the equations \eqref{3.1} hold for $\www b_1,\www b_2$,
they hold for any basis $b_1,b_2$ as in (a).

\medskip
(b) If $\kappa=1$ then \eqref{3.13} shows that $L$ is on
$Ml(f)_{1,\Z}$ up to the factor $\chi$ the standard symplectic
form. Therefore \eqref{3.2} holds.
If $\kappa<1$ then \eqref{3.13} shows \eqref{3.3}.

The restriction map $G_\Z\to \Aut(Ml(f)_{1,\Z})$ contains $T$.
This follows from \eqref{3.6} (whose proof below does not use 
this fact),
because obviously $(\delta,\alpha,\beta,\gamma)$ as in
\eqref{3.7} with $\delta=1$ exist.

This shows (b) in the case $\kappa<1$. In the case $\kappa=1$
we did not calculate lifts to $G_\Z$ of other elements of 
$\Aut(Ml(f)_{1,\Z},L)$. In this case the surjectivity of
$G_\Z\to \Aut(Ml(f)_{1,\Z},L)$ follows in two ways:
It follows from \cite[III.2.6]{Kl}, and it follows
from calculations in \cite{He1}, which are discussed
in the proof of theorem \ref{t6.1}.

\medskip
(c) We will prove (c) for the special choice $\www b_1,\www b_2$.
Then (c) holds for any $b_1,b_2$ as in (a) because by 
the surjectivity of the map $G_\Z\to\Aut(Ml(f)_{1,\Z},L)$,
an element $g\in G_\Z$ with $g(\www b_1)=b_1,g(\www b_2)=b_2$ exists.
Define
\begin{eqnarray}\label{3.14}
Ml_\Z^{[1]}&:=&\Z\www b_1\oplus\bigoplus_{i=1}^{p-1}\Z\delta_i,\ 
Ml_\C^{[1]}:=Ml_\Z^{[1]}\otimes_\Z\C,\\
Ml^{(1)}_\C&:=&Ml_\C^{[1]}\cap Ml(f)_{\neq 1},\label{3.15}
\end{eqnarray}
and analogously $Ml_\Z^{[2]},Ml_\C^{[2]},Ml^{(2)}_\C$ and 
$Ml_\Z^{[3]},Ml_\C^{[3]},Ml^{(3)}_\C$.

A look at the matrix $M_M$ shows the following.
\begin{eqnarray}\label{3.16}
M_h:\delta_1+\www b_1\mapsto \delta_2\mapsto ...\mapsto 
\delta_{p-1}\mapsto -(\delta_1+...+\delta_{p-1})\mapsto 
\delta_1+\www b_1.
\end{eqnarray}
Therefore $Ml_\Z^{[1]}$ is a cyclic $M_h$-module with
characteristic polynomial $t^p-1$, and $Ml^{[1]}_\C=\C\www b_1
\oplus Ml^{(1)}_\C$, and $M_h$ on $Ml^{(1)}_\C$ has the characteristic
polynomial $(t^p-1)/(t-1)$.

Lemma \ref{t2.5} applies and shows
\begin{eqnarray}\label{3.17}
\Aut(Ml_\Z^{[1]},L)=\{\pm(M_h|_{Ml^{[1]}_\Z})^\alpha\, |\, 
\alpha\in\{0,1,...,p-1\}\}.
\end{eqnarray}

Finally, $M_h,I$ and $L$ are well defined on the quotient lattice
$Ml_\Z^{[1]}/\Z\cdot\www b_1$, and 
$(Ml_\Z^{[1]}/\Z\cdot\www b_1,-I)$ is a root lattice of type 
$A_{p-1}$. The last statement follows immediately from the 
part of the Coxeter-Dynkin diagram which corresponds to 
$\delta_1,...,\delta_{p-1}$.

$Ml_\Z^{[2]}$ and $Ml_\Z^{[3]}$ have the same properties
as $Ml_\Z^{[1]}$, with $q$ respectively $r$ instead of $p$.

Now 
$$Ml(f)_{\neq 1}=Ml^{(1)}_\C\oplus Ml^{(2)}_\C\oplus Ml^{(3)}_\C$$
is clear. The $\Z$-lattice
$$Ml^{[1]}_\Z+Ml^{[2]}_\Z+Ml^{[3]}_\Z=\Z\www b_1\oplus
\bigoplus_{i=1}^{\mu-2}\Z\delta_i = 
(\C\www b_1\oplus Ml(f)_{\neq 1})\cap Ml(f)$$
is a primitive sublattice of $Ml(f)$ of rank $\mu-1$.
Any $g\in G_\Z$ with $g(\www b_1)=\pm \www b_1$ maps it
to itself, because it maps $\C\www b_1$ and $Ml(f)_{\neq 1}$
and $Ml(f)$ to themselves. $g$ maps also the quotient lattice
$$(Ml^{[1]}_\Z+Ml^{[2]}_\Z+Ml^{[3]}_\Z)/\Z\www b_1 = 
Ml^{[1]}_\Z/\Z\www b_1\oplus Ml^{[2]}_\Z/\Z\www b_1
\oplus Ml^{[3]}_\Z/\Z\www b_1$$
to itself.
But this is together with $-I$ an orthogonal sum of lattices
of types $A_{p-1},\ A_{q-1}$ and $A_{r-1}$.
Therefore $g$ can only permute the summands, and only those
summands of equal rank.

If $p=q$, a special element $\sigma_{12}\in G_\Z$ is given by
\begin{eqnarray*}
\sigma_{12}(\delta_i)&=& \delta_{p-1+i},\quad 
\sigma_{12}(\delta_{p-1+i})=\delta_i
\quad\textup{for }1\leq i\leq p-1,\\
\sigma_{12}(\delta_j)&=&\delta_j
\quad\textup{for }p+q-2\leq j\leq \mu.
\end{eqnarray*}
$\sigma_{12}\in G_\Z$ follows immediately from the symmetry
of the Coxeter-Dynkin diagram. 
Similarly $\sigma_{23}\in G_\Z$ is defined if $q=r$.
In any case, these elements generate a subgroup $U_2\subset G_\Z$
with the properties in (c).

Therefore, starting with an arbitrary element 
$\www g\in G_\Z$ if $\kappa<1$ respectively 
$\www g\in\{g\in G_\Z\, |\, g(\www b_1)=\pm\www b_1\}$ if $\kappa=1$,
one can compose it with $\pm\id$ and an element of $U_2$,
and one obtains an element $g\in G_\Z$ with $g(\www b_1)=\www b_1$
and $g(Ml^{[j]}_\Z)=Ml^{[j]}_\Z$ for $j=1,2,3$.
Then $g|_{Ml^{[1]}_\Z}=(M_h|_{Ml^{[1]}_\Z})^\alpha$ for a unique
$\alpha\in\{0,1,...,p-1\}$, and similarly
with $\beta\in\{0,1,...,q-1\}$ and $\gamma\in\{0,1,...,r-1\}$
for $Ml^{[2]}_\Z$ and $Ml^{[3]}_\Z$.
Also $g(\www b_2)=\www b_2+\delta\www b_1$ for some $\delta\in\Z$.
One calculates, observing \eqref{3.16},
\begin{eqnarray}\label{3.18}
M_h\left(\sum_{i=1}^{p-1}\frac{p-i}{p}\delta_i\right)
&=&\left(\sum_{i=1}^{p-1}\frac{p-i}{p}\delta_i\right)
-(\delta_1+\www b_1)+\frac{1}{p}\www b_1,\\
M_h^\alpha\left(\sum_{i=1}^{p-1}\frac{p-i}{p}\delta_i\right)
&=&\left(\sum_{i=1}^{p-1}\frac{p-i}{p}\delta_i\right)
-\left(\www b_1+\sum_{i=1}^\alpha\delta_i\right)
+\frac{\alpha}{p}\www b_1.\hspace*{1cm}\label{3.19}
\end{eqnarray}
The definition \eqref{3.11} of $\www b_2$ shows
\begin{eqnarray}\label{3.20}
-\delta_{\mu-1} =
-\frac{1}{\chi}\www b_2 + \sum_{i=1}^{p-1}\frac{p-i}{p}\delta_i
+\sum_{i=1}^{q-1}\frac{q-i}{q}\delta_{p-1+i}
+\sum_{i=1}^{r-1}\frac{r-i}{r}\delta_{p+q-2+i},
\end{eqnarray}
\eqref{3.19} gives
\begin{eqnarray}\label{3.21}
g(-\delta_{\mu-1})&=&-\delta_{\mu-1} + 
\left(\frac{-\delta}{\chi}+\frac{\alpha}{p}+\frac{\beta}{q}
+\frac{\gamma}{r}\right)\cdot\www b_1\\
&-&\left(\www b_1+\sum_{i=1}^\alpha\delta_i\right)
-\left(\www b_1+\sum_{i=p}^{p-1+\beta}\delta_i\right)
-\left(\www b_1+\sum_{i=p+q-1}^{p+q-2+\gamma}\delta_i\right) .
\nonumber
\end{eqnarray}
Therefore
\begin{eqnarray}\label{3.22}
\frac{\alpha}{p}+\frac{\beta}{q}+\frac{\gamma}{r}
&\equiv& \frac{\delta}{\chi}\mod 1\\
\textup{and}\quad g&=&
T^\delta\times (M_h|_{Ml^{(1)}_\C})^\alpha
\times (M_h|_{Ml^{(2)}_\C})^\beta
\times (M_h|_{Ml^{(3)}_\C})^\gamma. \nonumber
\end{eqnarray}
Thus $g\in U_1$, so $G_\Z\subset (U_1\rtimes U_2)\times\{\pm\id\}$.

Vice versa, we have to show $U_1\subset G_\Z$.
Fix a $g\in U_1$. It respects the decomposition
$$Ml(f)_\C = Ml(f)_1\oplus Ml^{(1)}_\C\oplus Ml^{(2)}_\C\oplus Ml^{(3)}_\C.$$
This is a left and right orthogonal decomposition with respect
to the Seifert form $L$.
The restriction of $g$ to each of the four blocks respects $L$ there,
so $g\in \Aut(Ml(f)_\C,L)$. It restricts on $Ml^{[1]}_\C$
to $M_h^\alpha$, so it maps the lattice $Ml^{[1]}_\Z$ to itself,
and analogously the lattices $Ml^{[2]}_\Z$ and $Ml^{[3]}_\Z$,
thus also the sum $Ml^{[1]}_\Z+Ml^{[2]}_\Z+Ml^{[3]}_\Z$.
This sum is a primitive sublattice of $Ml(f)$ of rank $\mu-1$
with 
$$Ml(f)=\left(Ml^{[1]}_\Z+Ml^{[2]}_\Z+Ml^{[3]}_\Z\right)\oplus
\Z\delta_{\mu-1}.$$
The calculation above of $g(-\delta_{\mu-1})$ shows
$g(\delta_{\mu-1})\in Ml(f)$ and 
$g(\delta_{\mu-1})\equiv\delta_{\mu-1}$ modulo the sublattice.
Therefore $g\in G_\Z$. \hfill$\Box$

\section{The group $G_\Z$ for 6 of the 28 exceptional
unimodal and bimodal singularities}\label{c4}

\noindent
The 14 1-parameter families of exceptional unimodal singularities
and the 14 2-parameter families of exceptional bimodal
singularities had been classified by Arnold.
Normal forms can be found in \cite{AGV1}.
In \cite[theorem 8.3]{He7} the group $G_\Z$ was calculated for 22
of the 28 families, namely for those families where
all eigenvalues of the monodromy have multiplicity 1.
In these cases it turned out that $G_\Z$ is simply
$\{\pm M_h^k\, |\, k\in\Z\}$. The proof used
lemma \ref{t2.5} and 
that the Milnor lattice is in these cases a cyclic
monodromy module. 

In this section $G_\Z$ will be determined
for the remaining 6 of the 28 families.
These are the families $Z_{12},Q_{12},U_{12},
Z_{18},Q_{16},U_{16}$.
In these cases, some eigenvalues have multiplicity 2.
This is similar to the case of the singularity $D_{2m}$,
which had also been treated in \cite[theorem 8.4]{He7}.
Also the proof will be similar. It will again use 
lemma \ref{t2.5} and combine that with an additional
analysis of the action of $G_\Z$ on the sum of the
eigenspaces with dimension $=2$, see lemma \ref{t4.2} below.

This lemma presents a surprise,
it points at a funny generalization of the
number theoretic fact
$\Z[e^{2\pi i/n}]\cap S_1=\{\pm e^{2\pi i k/n}\, |\, k\in\Z\}$.
Also, the lemma uses at the end a calculation of 
$G^{mar}\subset G_\Z$,
which will only come in section \ref{c7}.
The whole section \ref{c4} is devoted to the proof of the
following theorem.

\begin{theorem}\label{t4.1}
In the case of the 6 families of exceptional unimodal and bimodal
singularities $Z_{12},Q_{12},U_{12},Z_{18},Q_{16}$ and $U_{16}$,
the group $G_\Z$ is $G_\Z=\{\pm M_h^k\, |\, k\in\Z\}\times U$
with 
\begin{eqnarray}\label{4.1}
\begin{array}{l|l|l|l|l|l|l}
 & Z_{12} & Q_{12} & U_{12} & Z_{18} & Q_{16} & U_{16}\\ \hline
U\cong & \{\id\} & S_2 & S_3 & \{\id\} & S_2 & S_3
\end{array}
\end{eqnarray}
(This is independent of the number of variables, i.e.
it does not change under stabilization.)
\end{theorem}

{\bf Proof:}
Here we consider all 6 families as surface singularities.
Their characteristic polynomials $p_{ch}$ have
all the property $p_{ch}=p_1\cdot p_2$ with $p_2|p_1$
and $p_1$ having only simple roots.
They are as follows. Here $\Phi_m$ is the cyclotomic 
polynomial of primitive $m$-th unit roots.
\begin{eqnarray}\label{4.2}
\begin{array}{l|l|l|l|l|l|l}
 & Z_{12} & Q_{12} & U_{12} & Z_{18} & Q_{16} & U_{16}\\ \hline
p_{ch} & \Phi_{22}\Phi_2^2 &  \Phi_{15}\Phi_3^2 & 
\Phi_{12}\Phi_6\Phi_4^2\Phi_2^2 &\Phi_{34}\Phi_2^2 & 
\Phi_{21}\Phi_3^2 &  \Phi_{15}\Phi_5^2 \\
p_1 & \Phi_{22}\Phi_2 & \Phi_{15}\Phi_3 & 
\Phi_{12}\Phi_6\Phi_4\Phi_2 & \Phi_{34}\Phi_2 & \Phi_{21}\Phi_3 &
\Phi_{15}\Phi_5 \\
p_2 & \Phi_2 & \Phi_3 & \Phi_4\Phi_2 & \Phi_2 & \Phi_3 & \Phi_5
\end{array}
\end{eqnarray}

Orlik \cite{Or} had conjectured that the Milnor lattice of
any quasihomogeneous singularity is a sum of cyclic 
monodromy modules and that the characteristic polynomials
of $M_h$ on the cyclic pieces are $p_1,...,p_r$
where $p_{ch}=p_1\cdot ...\cdot p_r$ and $p_r|p_{r-1}|...|p_1$
and $p_1$ has simple roots ($r$ and $p_1,...,p_r$ are uniquely
determined by this). In the case of curve singularities,
Michel and Weber \cite{MW} have a proof of this conjecture.
In \cite[3.1]{He1} the conjecture 
was proved (using Coxeter-Dynkin diagrams)
for all those quasihomogeneous surface singularites of modality
$\leq 2$ which are not stabilizations of curve singularities.
So especially, the conjecture is true for the 
families of singularities 
$Z_{12},Q_{12},U_{12},Z_{18},Q_{16},U_{16}$.
There are $a_1,a_2\in Ml(f)$ with
\begin{eqnarray}\label{4.3}
Ml(f)=\left(\bigoplus_{i=0}^{\deg p_1-1}\Z\cdot M_h^i(a_1)\right)
\oplus \left(\bigoplus_{i=0}^{\deg p_2-1}\Z\cdot M_h^i(a_2)\right)
=: B_1\oplus B_2.
\end{eqnarray}

Denote
\begin{eqnarray}\label{4.4}
B_3:=\ker\left(p_2(M_h):Ml(f)_\C\to Ml(f)_\C\right)\cap Ml(f).
\end{eqnarray}
It is a primitive sublattice of $Ml(f)$ of rank $2\deg p_2$.
Also,
\begin{eqnarray}\label{4.5}
(B_1)_\C=\ker((p_1/p_2)(M_h))\oplus (B_1\cap B_3)_\C,\quad
B_2\subset B_3,
\end{eqnarray}
and $B_1\cap B_3$ and $B_2$ are both $M_h$-invariant 
primitive sublattices
of $B_3$ of rank $\deg p_2$. Together they generate $B_3$.

Any $g\in G_\Z$ with $g|_{B_3}=\pm (M_h|_{B_3})^k$ 
for some $k\in\Z$ 
restricts because of \eqref{4.5} to an automorphism of $B_1$.
Lemma \ref{t2.5} applies and shows
$g|_{B_1}=\pm (M_h|_{B_1})^l$ for some $l\in\Z$.
Now $g|_{B_3}=\pm (M_h|_{B_3})^k$ enforces 
$k\equiv l\mod \lcm(m\, |\, \Phi_m|p_2)$
and $g=\pm M_h^l$. Therefore
\begin{eqnarray}\label{4.6}
\{g\in G_\Z\, |\, g|_{B_3}=\pm (M_h|_{B_3})^k\textup{ for some }k\}
=\{\pm M_h^k\, |\, k\in\Z\}.
\end{eqnarray}

Lemma \ref{t4.2} (c) determines $\Aut(B_3,L)$,
\begin{eqnarray}\label{4.7}
\Aut(B_3,L)=\{\pm (M_h|_{B_3})^k\, |\, k\in\Z\}
\times U\quad\textup{ with }U\textup{ as in }\eqref{4.1}.
\end{eqnarray}
Lemma \ref{t4.2} (d) shows that the map 
$G_\Z\to \Aut(B_3,L)$ is surjective.
Together with \eqref{4.6} this gives \eqref{4.1}.
\hfill$\Box$

\begin{lemma}\label{t4.2}
(a) Let $V_\Z$ be a $\Z$-lattice of rank $2$ with a $\Z$-basis
$\uuuu{b}=(b_1,b_2)$ and a symmetric pairing $L_\Z$ given by
$$L_\Z(\uuuu{b}^t,\uuuu{b})=
\begin{pmatrix}2&-1\\-1&m\end{pmatrix}\quad\textup{for some }
m\in\N_{\geq 2}.$$
Define $\xi:=e^{2\pi i/l}$ for some $l\in\{3,4,5\}$,
where we exclude the cases $(m\geq 3,l=5)$. 
Define $V_\C:=V_\Z\otimes_\Z\C$,
$V_{\Z[\xi]}:=V_\Z\otimes_\Z\Z[\xi]\subset V_\C$,
and extend $L_\Z$ sesquilinearly (=linear$\times$semilinear)
to $V_\C$. 
Then
\begin{eqnarray}\label{4.8}
&&\{r\in V_{\Z[\xi]}\, |\, L_\C(r,r)=2\}\\
&=& \{\pm \xi^k\, |\, k\in\Z\}\times 
\{r\in V_\Z\, |\, L_\Z(r,r)=2\}.\nonumber
\end{eqnarray}
In the case $m\geq 3$
\begin{eqnarray}\label{4.9}
&&\{r\in V_{\Z[\xi]}\, |\, L_\C(r,r)=m,r\notin\Z[\xi]b_1\}\\
&=& \{\pm \xi^k\, |\, k\in\Z\}\times 
\{r\in V_\Z\, |\, L_\Z(r,r)=m,r\notin\Z b_1\}.\nonumber
\end{eqnarray}

\medskip
(b) In the situation of (a)
\begin{eqnarray}\label{4.10}
\Aut(V_{\Z[\xi]},L_\C)&=&\{\pm\xi^k\, |\, k\in\Z\}\cdot
\Aut(V_\Z,L_\Z)\\ 
&&(\pm\id\textup{ lives on both sides, therefore not }\times)
\nonumber\\ 
\nonumber
\Aut(V_\Z,L_\Z)&\cong&
\{\pm\begin{pmatrix}1&0\\0&1\end{pmatrix},
\pm\begin{pmatrix}-1&1\\0&1\end{pmatrix}\}\\
&\cong& \{\pm\id\}\times S_2\quad\textup{in the cases }m\geq 3,
\label{4.11}\\
\Aut(V_\Z,L_\Z)&\cong&
\Aut(\textup{root lattice of type }A_2)\nonumber \\
&\cong& \{\pm\id\}\times S_3\quad\textup{in the case }m=2.
\label{4.12}
\end{eqnarray}

\medskip
(c) In the situation of the proof of theorem \ref{t4.1}
\begin{eqnarray}\label{4.13}
\Aut(B_3,L)=\{\pm (M_h|_{B_3})^k\, |\, k\in\Z\}
\times U\quad\textup{ with }U\textup{ as in }\eqref{4.1}.
\end{eqnarray}

\medskip
(d) In the situation of the proof of theorem \ref{t4.1}
the map $G_\Z\to\Aut(B_3,L)$ is surjective.
\end{lemma}

{\bf Proof:}
(a) An element 
$r=r_1b_1+r_2b_2$ with $r_1,r_2\in\Z[\xi]$ satisfies
\begin{eqnarray}
L_\C(r,r)&=&2|r_1|^2-(r_1\oooo{r_2}+\oooo{r_1}r_2)+m|r_2|^2
\nonumber \\
&=& |r_1|^2+|r_1-r_2|^2+(m-1)|r_2|^2.\label{4.14}
\end{eqnarray}
First consider the cases $l\in\{3,4\}$. Then 
$\Z[\xi]\cap \R=\Z$. Then the three absolute values in \eqref{4.14}
are non-negative integers. Their sum is 2 if and only if
\begin{eqnarray}\label{4.15}
|r_1|=1,r_2=0&&\textup{ in the cases }m\geq 3,\\
\left. \begin{array}{l}
(|r_1|,|r_2|)\in\{(1,0),(0,1),(1,1)\}\\
\textup{and in the last case }r_1=r_2
\end{array}\right\}
&&\textup{ in the case }m=2.\label{4.16}
\end{eqnarray}
In the case $m\geq 3$ and in the case of an $r\notin\Z[\xi]b_1$,
the sum of the three absolute values in \eqref{4.15}
is $m$ if and only if
\begin{eqnarray}\label{4.17}
(r_1=0,|r_2|=1)\textup{ or } (r_1=r_2,|r_1|=1).
\end{eqnarray}
Together with $\Z[\xi]\cap S^1=\{\pm \xi^k\, |\, k\in\Z\}$
this shows part (a) in the cases $l\in\{3,4\}$.

It rests to consider the case $(m,l)=(2,5)$. 
In that case write
$$r_1=r_{10}+r_{11}\xi+r_{12}\xi^2+r_{13}\xi^3,\quad
r_2=r_{20}+r_{21}\xi+r_{22}\xi^2+r_{23}\xi^3$$
with $r_{ij}\in\Z$.
Then
\begin{eqnarray*}
L_\C(r,r)&=&
2|r_1|^2+2|r_2|^2-(r_1\oooo{r_2}+\oooo{r_1}r_2)\\
&=& 2\left[
\sum_{j=0}^3r_{1j}^2 + (\xi+\xi^4)\sum_{j=1}^3r_{1j}r_{1,j-1}
+(\xi^2+\xi^3)\sum_{j-k\geq 2}r_{1j}r_{1k}\right]
\\
&+& 2\left[
\sum_{j=0}^3r_{2j}^2 + (\xi+\xi^4)\sum_{j=1}^3r_{2j}r_{2,j-1}
+(\xi^2+\xi^3)\sum_{j-k\geq 2}r_{2j}r_{2k}\right]
\\
&-&\left[2\sum_{j=0}^3r_{1j}r_{2j}
+(\xi+\xi^4)\sum_{j=1}^3(r_{1j}r_{2,j-1}+r_{1,j-1}r_{2j}) \right.
\\
&& + \left. 
(\xi^2+\xi^3)\sum_{j-k\geq 2}(r_{1j}r_{2k}+r_{1k}r_{2j})\right]\\
&=& A_1+A_2\cdot\frac{\sqrt{5}}{2}
\quad\textup{with }A_1,A_2\in\Z.
\end{eqnarray*}
It is not so easy to find, but easy to check that $A_1$ is equal to
\begin{eqnarray}
&&\frac{1}{4}\sum_{j=0}^3\left[r_{1j}^2+r_{2j}^2
+(r_{1j}-r_{2j})^2\right] \label{4.18}\\
&+&\frac{1}{4}\sum_{j<k}\left[
(r_{1j}-r_{1k}-r_{2j}+r_{2k})^2
+(r_{1j}-r_{1k})^2 
+(r_{2j}-r_{2k})^2\right]. \nonumber
\end{eqnarray}
Now it is an easy exercise to find the 8-tuples
$(r_{10},...,r_{23})\in\Z^8$ for which \eqref{4.18} 
takes the value 2. They are (here $e_j=(\delta_{ij})_{i=1,...,8}$
for $j=1,...,8$ is the standard basis of $\Z^8$) 
\begin{eqnarray}\label{4.19}
\pm e_1,...,\pm e_8,\pm (e_1+e_5),\pm (e_2+e_6),\pm (e_3+e_7),
\pm (e_4+e_8),\\
\pm (1,1,1,1,0,0,0,0),\pm (0,0,0,0,1,1,1,1),\pm (1,1,1,1,1,1,1,1).
\nonumber
\end{eqnarray}
Observe $1+\xi+\xi^2+\xi^3=-\xi^4$ and
$$\{r\in V_\Z\, |\, L_\Z(r,r)=2\}=\{\pm b_1,\pm b_2,\pm(b_1+b_2)\}.$$
The coefficients $(r_{10},...,r_{23})$ in \eqref{4.19}
give precisely the elements $r=r_1b_1+r_2b_2$ on the
right hand side of \eqref{4.8}. This shows part (a) for
$(m,l)=(2,5)$.

\medskip
(b) Any element $g$ of $\Aut(V_{\Z[\xi]},L_\C)$ will map
the sets in \eqref{4.8} and \eqref{4.9} to themselves. 
The basis elements 
$b_1$ and $b_2$ are mapped to two elements in these sets 
with $L_\C(g(b_1)),g(b_2))=L_\Z(b_1,b_2)=-1$. 
Therefore $g$ is up to 
a factor in $\{\pm\xi^k\, |\, k\in\Z\}$ an element
of $\Aut(V_\Z,L_\Z)$. This shows \eqref{4.10}.

In the case $m\geq 3$ 
\begin{eqnarray*}
\{r\in V_\Z\, |\, L_\Z(r,r)=2\}&=&\{\pm b_1\}\\
\textup{and}\quad
\{r\in V_\Z\, |\, L_\Z(r,r)=m\}&=&\{\pm b_2,\pm(b_1+b_2)\}.
\end{eqnarray*}
This shows \eqref{4.11}. The case $m=2$ is the case of
the root lattice of type $A_2$. \eqref{4.12} is well known
and easy to see.

\medskip
(c) In the cases $Z_{12}$ and $Z_{18}$ as curve singularities,
the examples \ref{t2.2} (ii) and (iii) showed
$\Aut(\Rad(I),L)=\{\pm\id\}$. Here 
$\Rad(I)=\ker(M_h^{curve\ case}-\id)$. Under stabilization
$\Rad(I)$ becomes $B_3$, and $L$ changes just the sign,
see \eqref{2.16}. Thus $\Aut(B_3,L)=\{\pm\id\}$.
Because of $M_h|_{B_3}=-\id$ and $U:=\{\id\}$ in the cases
$Z_{12}$ and $Z_{18}$, this shows (c) in these cases.

Now consider the cases $Q_{12},Q_{16},U_{12}$ and $U_{16}$.
Here part (b) will be used, but that has to be prepared.

The normal forms of the quasihomogeneous surface singularities, 
given in section \ref{c7}, show that they are sums 
of singularities in different variables 
of types $A_l$ and $D_{2m}$ with $(l,2m)$ as follows,
\begin{eqnarray}\label{4.20}
\begin{array}{l|l|l|l|l}
 &Q_{12}&Q_{16}&U_{12}&U_{16}\\ \hline
(l,2m)&(2,6)&(2,8)&(3,4)&(4,4)\\
 & A_2\otimes D_6 & A_2\otimes D_8 &A_3\otimes D_4 & A_4\otimes D_4
\end{array}
\end{eqnarray}
Here the singularity $A_l$ is in one variable
and has the characteristic polynomial 
$p_{ch}^{A_l}=(t^{l+1}-1)/(t-1)$,
the singularity $D_{2m}$ is a curve singularity
and has the characteristic polynomial 
$p_{ch}^{D_{2m}}=(t^{2m-1}-1)\Phi_1$.
The Thom-Sebastiani results which were cited in section \ref{c2}
apply,
\begin{eqnarray}\label{4.21}
(Ml(f),L)&\cong& (Ml(A_l),L_{A_l})\otimes (Ml(D_{2m}),L_{D_{2m}}),\\
M_h &\cong& M_h^{A_l}\otimes M_h^{D_{2m}},\nonumber
\end{eqnarray}
and show
\begin{eqnarray}
p_2 &=& p_{ch}^{A_l},\nonumber \\
(B_3,L) &\cong & (Ml(A_l),L_{A_l})\otimes 
(Ml(D_{2m})_{1,\Z},L_{D_{2m}}),
\label{4.22}\\
M_h|_{B_3} &\cong & M_h^{A_l}\otimes \id.\nonumber
\end{eqnarray}
The pair $(Ml(D_{2m})_{1,\Z},L_{D_{2m}})$
was considered in example \ref{t2.2} (i).
There is a $\Z$-basis $\uuuu{b}=(b_1,b_2)$ of 
$Ml(D_{2m})_{1,\Z}$ with
\begin{eqnarray}\label{4.23}
L_{D_{2m}}(\uuuu{b}^t,\uuuu{b})=
\begin{pmatrix}-2&1\\1&-m\end{pmatrix}.
\end{eqnarray}
The pairings $L$ and $L_{A_l}\otimes L_{D_{2m}}$
will be extended sesquilinearly from the $\Z$-lattices
to the $\C$-vector spaces.

The $\Z$-lattice $Ml(A_{l})$ is a cyclic monodromy module.
Choose a generator $e$ of it. Therefore 
$Ml(A_l)\otimes Ml(D_{2m})_{1,\Z}$ is a sum of two 
cyclic monodromy modules, and generators are
$e\otimes b_1$ and $e\otimes b_2$.
For any automorphism $g$ of 
$(Ml(A_l)\otimes Ml(D_{2m})_{1,\Z},M_h^{A_l}\otimes\id)$ 
there are unique polynomials
$g_1,g_2,g_3,g_4\in\Z[t]$ of degree $\leq \deg p_2-1$
such that
\begin{eqnarray}\label{4.24}
\begin{pmatrix}g(v\otimes b_1)\\g(v\otimes b_2)\end{pmatrix}
=\begin{pmatrix}g_1(M_h^{A_l})(v)\otimes b_1 + 
g_3(M_h^{A_l}) (v)\otimes b_2\\
g_2(M_h^{A_l})(v)\otimes b_1 + 
g_4(M_h^{A_l}) (v)\otimes b_2\end{pmatrix}
\end{eqnarray}
for any $v\in Ml(A_l)$.

Now choose any eigenvalue $\xi$ of $M_h^{A_l}$.
Then $\Z[\xi]$ is a principal ideal domain.
The space $\ker(M_h^{A_l}-\xi\id)\cap Ml(A_l)_{\Z[\xi]}$
is a free $\Z[\xi]$-module of rank 1. 
Choose a generating vector $v$. This choice gives an 
isomorphism from this space to $\Z[\xi]$.
The spaces
$$\ker(M_h-\xi\id)\cap Ml(f)_{\Z[\xi]}
\cong (\ker(M_h^{A_l}-\xi\id)\cap Ml(A_l)_{\Z[\xi]})
\otimes Ml(D_{2m})_{1,\Z[\xi]}$$
are free $\Z[\xi]$-modules of rank 2.
The space on the right hand side has the $\Z[\xi]$-basis
$(v\otimes b_1,v\otimes b_2)=:v\otimes\uuuu{b}$.
Now \eqref{4.24} becomes
\begin{eqnarray}\label{4.25}
g(v\otimes\uuuu{b})=v\otimes \uuuu{b}\cdot
\begin{pmatrix}g_1(\xi)&g_2(\xi)\\
g_3(\xi)&g_4(\xi)\end{pmatrix}.
\end{eqnarray}
The pairing satisfies 
\begin{eqnarray}\label{4.26}
(L_{A_l}\otimes L_{D_{2m}})(v\otimes \uuuu{b})
= L_{A_l}(v,v)\cdot 
\begin{pmatrix}-2&1\\1&-m\end{pmatrix},
\end{eqnarray}
where $L_{A_l}(v,v)\in\Z[\xi]\cap\R_{>0}$.
This space with this pairing is up to a scalar isomorphic
to a pair $(V_{\Z[\xi]},L_\C)$ considered in the parts (a) 
and (b).
Therefore by part (b), its group of automorphisms is isomorphic
to $\{\pm\xi^k\, |\, k\in\Z\}\cdot \Aut(V_\Z,L_\Z)$.

Thus any element of $\Aut(B_3,L)$ restricts on
$\ker(M_h-\xi\id)\cap Ml(f)_{\Z[\xi]}$ to such an automorphism.
In the cases $Q_{12},Q_{16}$ and $U_{16}$ the polynomial
$p_2=\Phi_3,\Phi_3$ respectively $\Phi_5$ is irreducible, 
so all its zeros $\xi$ are conjugate. 
Therefore then 
\begin{eqnarray*}
\Aut(B_3,L)
&\cong& \{\pm(M_h|_{B_3})^k\, |\, k\in\Z\} \cdot\Aut(V_\Z,L_\Z\}\\
&\cong& \{\pm(M_h|_{B_3})^k\, |\, k\in\Z\} \times U,
\end{eqnarray*}
which proves (c) in these cases.

In the case $U_{12}$ the characteristic polynomial 
$p_{ch}^{A_3}=p_2=\Phi_4\Phi_2$ is reducible. 
Consider an automorphism $g$ of 
$$(Ml(A_l)\otimes Ml(D_{2m})_{1,\Z},L_{A_l}\otimes L_{D_{2m}}).$$
It is determined by the polynomials $g_1,g_2,g_3,g_4\in\Z[t]$
in \eqref{4.24}. 
For $\xi=i$ and for $\xi=-1$ it gives an automorphism
of $\Z[\xi]v\otimes b_1\oplus\Z[\xi]v\otimes b_2$
which is given by a matrix 
$\begin{pmatrix}g_1(\xi)&g_2(\xi)\\ g_3(\xi)&g_4(\xi)\end{pmatrix}$
which is by part (b) in 
\begin{eqnarray}\label{4.27}
\{\pm\xi^k\, |\, k\in\Z\}&\cdot& \left\{
\pm\begin{pmatrix} 1&0 \\ 0&1 \end{pmatrix},
\pm\begin{pmatrix} 0&-1 \\ 1&-1 \end{pmatrix},
\pm\begin{pmatrix} -1&1 \\ -1&0 \end{pmatrix}\right. ,\\
&&\left. \pm\begin{pmatrix} 0&1 \\ 1&0 \end{pmatrix},
\pm\begin{pmatrix} -1&0 \\ -1&1 \end{pmatrix},
\pm\begin{pmatrix} 1&-1 \\ 0&-1 \end{pmatrix}
\right\}.\nonumber
\end{eqnarray}
By multiplying $g$ with a suitable automorphism we can suppose
that the matrix for $\xi=i$ is the identity matrix.
Then 
\begin{eqnarray*}
\begin{pmatrix}g_1&g_2\\ g_3&g_4\end{pmatrix}
=\begin{pmatrix}1+(t^2+1)\www g_1&(t^2+1)\www g_2\\ 
(t^2+1)\www g_3&1+(t^2+1)\www g_4\end{pmatrix},
\end{eqnarray*}
for some $\www g_1,\www g_2,\www g_3,\www g_4\in\Z[t]$, so
\begin{eqnarray*}
\begin{pmatrix}g_1(-1)&g_2(-1)\\ g_3(-1)&g_4(-1)\end{pmatrix}
=\begin{pmatrix}1+2\www g_1(-1)&2\www g_2(-1)\\ 
2\www g_3(-1)&1+2\www g_4(-1)\end{pmatrix}.
\end{eqnarray*}
The only two possibilities are 
$\pm\begin{pmatrix}1&0\\0&1\end{pmatrix}$.
In the case of a minus sign 
$g\circ ((M_h^{A_l})^2\otimes \id_{D_{2m}})=-\id$,
in the case of a plus sign $g=\id$. This finishes
the proof of (c) in the case $U_{12}$.

\medskip
(d) In section \ref{c7} a subgroup $G^{mar}$ of $G_\Z$ 
will be calculated, and it will be shown that
$G^{mar}=\{\pm M_h^k\, |\, k\in\Z\}\times U$.
This shows $G_\Z\supset \{\pm M_h^k\, |\, k\in\Z\}\times U$
and that the map $G_\Z\to\Aut(B_3,L)$ is surjectiv.
\hfill$\Box$

\begin{remark}\label{t4.3}
The number theoretic fact
$\Z[e^{2\pi i/a}]\cap S^1=\{\pm e^{2\pi i k/a}\, |\, k\in\Z\}$
can be interpreted as saying that in the case of the 
$A_1$-lattice $V_\Z=\Z$ with $\Z$-basis $b_1=1$ and standard bilinear
form $L_\Z$ with $L_\Z(b_1,b_1)=1$ and hermitian extension
$L_\C$ to $\C$ the analogue of \eqref{4.8}
holds. Now \eqref{4.8} for $m=2$ 
can be seen as a generalization from the
case $A_1$ to the case $A_2$. Above it is proved only in the
cases $l=3,4,5$.
\end{remark}

\section{Review on $\mu$-constant monodromy groups $G^{mar}$,
marked singularities, their moduli spaces $M_\mu^{mar}$,
and Torelli type conjectures}\label{c5}

\noindent
This paper is a sequel to \cite{He7}.
That paper studied holomorphic functions germs
$f:(\C^{n+1},0)\to (\C,0)$ with an isolated singularity
at 0 from a global perspective.
Here we review most of the notions and results from \cite{He7}.

It defined the notions of {\it marked singularity}
and  {\it strongly marked singularity}. 
The marking uses the Milnor lattice $Ml(f)\cong \Z^\mu$
and the Seifert form $L$ on it, which are explained in
section \ref{c2}.

\begin{definition}\label{t5.1}
Fix one reference singularity $f_0$.

(a)
Then a strong marking for any singularity $f$
in the $\mu$-homotopy class of $f_0$ 
(i.e. there is a 1-parameter family of singularities with
constant Milnor number connecting $f$ and $f_0$)
is an isomorphism $\rho:(Ml(f),L)\to (Ml(f_0),L)$.

(b)
The pair $(f,\rho)$ is a {\it strongly marked singularity}.
Two strongly marked singularities $(f_1,\rho_1)$ and
$(f_2,\rho_2)$ are right equivalent (notation: $\sim_\RR$) 
if a coordinate change
$\varphi:(\C^{n+1},0)\to(\C^{n+1},0)$ with
$$f_1=f_2\circ \varphi\quad\textup{and}\quad
\rho_1=\rho_2\circ\varphi_{hom}$$
exists, where $\varphi_{hom}:(Ml(f_1),L)\to (Ml(f_2),L)$ 
is the induced isomorphism.

(c)
The notion of a marked singularity is slightly weaker.
If $f$ and $\rho$ are as above, then the pair ($f,\pm \rho)$
is a {\it marked singularity} (writing $\pm\rho$, the set
$\{\rho,-\rho\}$ is meant, neither $\rho$ nor 
$-\rho$ is preferred). 

(d)
Two marked singularities $(f_1,\rho_1)$ and $(f_2,\rho_2)$ 
are right equivalent (notation: $\sim_\RR$) 
if a coordinate change $\varphi$ with 
$$f_1=f_2\circ \varphi\quad
\textup{and}\quad\rho_1=\pm\rho_2\circ\varphi_{hom}$$
exists.
\end{definition}

\begin{remarks}\label{t5.2}
(i) The notion of a marked singularity behaves better than the
notion of a strongly marked singularity, because it is not
known whether all $\mu$-homotopy families of singularities
satisfy one of the following two properties:
\begin{eqnarray}\label{5.1}
\textup{Assumption (5.1):}&&\textup{Any singularity in the }
\mu\textup{-homotopy}\\
&&\textup{class of }f_0\textup{ has multiplicity }\geq 3.\nonumber\\
\textup{Assumption (5.2):}&&\textup{Any singularity in the }
\mu\textup{-homotopy} \label{5.2}\\
&&\textup{class of }f_0\textup{ has multiplicity }2.\nonumber
\end{eqnarray}
We expect that always one of two assumptions holds.
For curve singularities and singularities right equivalent
to semiquasihomogeneous singularities this is true, but
in general it is not known. In a $\mu$-homotopy family
where neither of the two assumptions holds, strong marking
behaves badly, see (ii).

\medskip
(ii) If $\textup{mult}(f)=2$ then $(f,\rho)\sim_\RR (f,-\rho)$,
which is easy to see. If $\textup{mult}(f)\geq 3$, then
$(f,\rho)\not\sim_\RR(f,-\rho)$, whose proof in \cite{He7}
is quite intricate. These properties imply that the
moduli space for strongly marked singularities discussed below
is not Hausdorff in the case of a $\mu$-homotopy class 
which satisfies neither one of the assumptions \eqref{5.1} 
or \eqref{5.2}
\end{remarks}

In \cite{He6} for the $\mu$-homotopy class of any singularity $f_0$
a moduli space $M_\mu(f_0)$ was constructed. As a set it is simply
the set of right equivalence classes of singularities in the
$\mu$-homotopy class of $f_0$. But in \cite{He6}
it is constructed as an analytic geometric quotient,
and it is shown that it is locally isomorphic to the 
$\mu$-constant stratum of a singularity modulo the action of
a finite group. The $\mu$-constant stratum of a singularity
is the germ $(S_\mu,0)\subset (M,0)$ within the germ of the
base space of a universal unfolding $F$ of $f$, such that for
a suitable representative 
\begin{eqnarray}\label{5.3}
S_\mu=\{t\in M\, |\, F_t\textup{ has only one singularity }x_0
\textup{ and }F_t(x_0)=0\}.
\end{eqnarray}
It comes equipped with a canonical complex structure,
and $M_\mu$ inherits a canonical structure, see the chapters
12 and 13 in \cite{He6}.

In \cite{He7} analogous results for marked singularities 
were proved. A better property is that $M^{mar}_\mu$ is
locally isomorphic to a $\mu$-constant stratum without
dividing out a finite group action. 
Therefore one can consider it as a {\it global $\mu$-constant
stratum} or as a {\it Teichm\"uller space for singularities}.
The following theorem
collects results from \cite[theorem 4.3]{He7}.

\begin{theorem}\label{t5.3}
Fix one reference singularity $f_0$. Define the sets
\begin{eqnarray}\label{5.4}
M^{smar}_\mu(f_0) &:=& \{\textup{strongly marked }(f,\rho)\, |\, \\
&&f\textup{ in the }\mu\textup{-homotopy class of }f_0\}/\sim_\RR,
\nonumber \\
M^{mar}_\mu(f_0) &:=& \{\textup{marked }(f,\pm\rho)\, |\,  
\label{5.5} \\
&&f\textup{ in the }\mu\textup{-homotopy class of }f_0\}/\sim_\RR.
\nonumber
\end{eqnarray}

(a) $M^{mar}_\mu(f_0)$ carries a natural canonical complex structure.
It can be constructed with the underlying reduced complex
structure as an analytic geometric quotient 
(see \cite[theorem 4.3]{He7} for details).

(b) The germ $(M^{mar}_\mu(f_0),[(f,\pm\rho)])$ with its canonical
complex structure is isomorphic to the $\mu$-constant stratum 
of $f$ with its canonical complex structure 
(see \cite[chapter 12]{He6} for the definition of that).

(c) For any $\psi\in G_\Z(f_0)=:G_\Z$, the map
$$\psi_{mar}:M_\mu^{mar}\to M_\mu^{mar},\quad
[(f,\pm\rho)]\to [(f,\pm\psi\circ\rho)]$$
is an automorphism of $M_\mu^{mar}$.  The action
$$G_\Z\times M_\mu^{mar}\to M_\mu^{mar},\quad 
(\psi,[(f,\pm\rho)]\mapsto \psi_{mar}([(f,\pm\rho)])$$
is a group action from the left.

(d) The action of $G_\Z$ on $M_\mu^{mar}$ is properly discontinuous.
The quotient $M_\mu^{mar}/G_\Z$ is the moduli space $M_\mu$ 
for right equivalence classes in the 
$\mu$-homotopy class of $f_0$, with its canonical complex structure.
Especially, $[(f_1,\pm\rho_1)]$ and $[(f_2,\pm\rho_2)]$ are in one
$G_\Z$-orbit if and only if $f_1$ and $f_2$ are right equivalent.

(e) If assumption \eqref{5.1} or \eqref{5.2} holds then 
(a) to (d) 
are also true for $M_\mu^{smar}$ and $\psi_{smar}$
with $\psi_{smar}([(f,\rho)]):=[(f,\psi\circ\rho)]$.
If neither\eqref{5.1} nor\eqref{5.2} holds then the natural
topology on $M^{smar}_\mu$ is not Hausdorff.
\end{theorem}

We stick to the situation in theorem \ref{t5.3} and
define two subgroups of $G_\Z(f_0)$.
The definitions in 
\cite[definition 3.1]{He7} are different, they use 
$\mu$-constant families. The following definitions are
a part of theorem 4.4 in \cite{He7}.

\begin{definition}\label{t5.4}
Let $(M^{mar}_\mu)^0$ be the topological component
of $M^{mar}_\mu$ (with its reduced complex structure) which
contains $[(f_0,\pm\id)]$. Then
\begin{eqnarray}\label{5.6}
G^{mar}(f_0)&:=& \{\psi\in G_\Z\, |\, \psi\textup{ maps }
(M_\mu^{mar})^0 \textup{ to itself}\}\subset G_\Z(f_0).
\end{eqnarray}
If assumption \eqref{5.1} or \eqref{5.2} holds, 
$(M^{smar}_\mu)^0$ and 
$G^{smar}(f_0)\subset G_\Z(f_0)$ are defined analogously.
\end{definition}

The following theorem is also proved in \cite{He7}.

\begin{theorem}\label{t5.5}
(a) In the situation above the map
\begin{eqnarray*}
G_\Z/G^{mar}(f_0)&\to& \{\textup{topological components of }
M_\mu^{mar}\}\\
\psi\cdot G^{mar}(f_0)&\mapsto& \textup{the component }
\psi_{mar}((M_\mu^{mar})^0)
\end{eqnarray*}
is a bijection.

(b) If assumption \eqref{5.1} or \eqref{5.2} holds then (a) 
is also true for $M_\mu^{smar}$ and $G^{smar}(f_0)$.

(c) $-\id\in G_\Z$ acts trivially on $M_\mu^{mar}(f_0)$. 
Suppose additionally that assumption \eqref{5.1} holds for $f_0$. 
Then $\{\pm\id\}$ acts freely on $M_\mu^{smar}(f_0)$,
and the quotient map
$$M_\mu^{smar}(f_0) \stackrel{/\{\pm\id\}}{\longrightarrow}
M_\mu^{mar}(f_0),\quad [(f,\rho)]\mapsto [(f,\pm\rho)]$$
is a double covering.
\end{theorem}

The first main conjecture in \cite{He7} is part (a) of the following
conjecture (the second main conjecture in \cite{He7} is conjecture 
\ref{t5.11} (a) below).

\begin{conjecture}\label{t5.6}
(a) Fix a singularity $f_0$. Then $M^{mar}_\mu$ is
connected. Equivalently (in view of theorem \ref{t5.5} (a)):
$G^{mar}(f_0)=G_\Z.$

(b) If the $\mu$-homotopy class of $f_0$ satisfies 
assumption \eqref{5.1}, then $-\id\notin G^{smar}(f_0)$.
\end{conjecture}

If (a) holds then (b) is equivalent to $M^{smar}_\mu(f_0)$
having two components.
If (a) and (b) hold, there should be a natural invariant which
distinguishes the index 2 subgroup $G^{smar}\subset G^{mar}=G_\Z$.
Anyway, part (a) is the more important conjecture.
Using the other definition of $G^{mar}$ in \cite{He7},
it says that up to $\pm\id$, any element of $G_\Z$
can be realized as  transversal monodromy of a
$\mu$-constant family with parameter space $S^1$.

The whole conjecture \ref{t5.6} had been proved 
in \cite{He7} for the simple singularities and those 22
of the 28 exceptional unimodal and bimodal singularities,
where all eigenvalues of the monodromy have only multiplicity
one \cite[theorems 8.3 and 8.4]{He7}. In this paper it
will be proved for the remaining unimodal and exceptional
bimodal singularities.

In order to understand the stabilizers
$\textup{Stab}_{G_\Z}([(f,\rho)])$ and 
$\textup{Stab}_{G_\Z}([(f,\pm\rho)])$
of points $[(f,\rho)]\in M^{smar}_\mu(f_0)$ and 
$[(f,\pm\rho)]\in M^{mar}_\mu(f_0)$,
we have to look at
the {\it symmetries} of a single singularity.
These had been discussed in \cite[chapter 13.2]{He6}.
The discussion had been taken up again in \cite{He7}.

\begin{definition}\label{t5.7}
Let $f_0=f_0(x_0,...,x_n)$ be a reference singularity 
and let $f$ be any singularity in the $\mu$-homotopy class of $f_0$.
If $\rho$ is a marking, then $G_\Z(f)=\rho^{-1}\circ G_\Z\circ\rho$.

We define
\begin{eqnarray}\label{5.7}
\RR &:=&\{\varphi:(\C^{n+1},0)\to (\C^{n+1},0)\quad
\textup{biholomorphic}\},\\
\RR^f &:=& \{\varphi\in\RR\, |\, f\circ \varphi=f\},\label{5.8}\\
R_f&:=& j_1\RR^f/(j_1\RR^f)^0,\label{5.9}\\
G^{smar}_\RR(f)&:=&\{\varphi_{hom}\, |\, \varphi\in\RR^f\}
\subset G_\Z(f),\label{5.10}\\
G^{mar}_\RR(f)&:=& \{\pm\psi\, |\, \psi\in G^{smar}_\RR(f)\}.
\label{5.11}
\end{eqnarray}
\end{definition}

Again, the definition of $G^{smar}_\RR$ is different from
the definition in \cite[definition 3.1]{He7}.
The characterization in \eqref{5.10} is 
\cite[theorem 3.3. (e)]{He7}.
$R_f$ is the finite group of components of the group
$j_1\RR^f$ of 1-jets of coordinate changes which leave $f$
invariant.
The following theorem collects results from several theorems
in \cite{He7}.

\begin{theorem}\label{t5.8}
Consider the data in definition \ref{t5.7}.

(a) If $\textup{mult}(f)\geq 3$ then $j_1\RR^f= R_f$.

(b) The homomorphism $()_{hom}:\RR^f\to G_\Z(f)$
factors through $R_f$. Its image is $(R_f)_{hom}=G^{smar}_\RR(f)
\subset G_\Z(f)$.

(c) The homomorphism $()_{hom}:R_f\to G^{smar}_\RR(f)$ 
is an isomorphism.

(d) 
\begin{eqnarray}\label{5.12}
-\id\notin G^{smar}_\RR(f)\iff \mult f\geq 3.
\end{eqnarray}
Equivalently: $G^{mar}_\RR(f)=G^{smar}_\RR(f)$ if $\mult f=2$, and
$G^{mar}_\RR(f)=G^{smar}_\RR(f)\times\{\pm\id\}$ if $\mult f\geq 3$.

(e) $G^{mar}_\RR(f)=G^{mar}_\RR(f+x_{n+1}^2)$.

(f) $M_h\in G^{smar}(f)$. If $f$ is quasihomogeneous then 
$M_h\in G^{smar}_\RR(f)$.

(g) For any $[(f,\rho)]\in M_\mu^{smar}$
\begin{eqnarray}\label{5.13}
\Stab_{G_\Z}([(f,\rho)]) &=& \rho\circ G^{smar}_\RR(f)\circ 
\rho^{-1},\\
\Stab_{G_\Z}([(f,\pm\rho)]) 
&=& \rho\circ G^{mar}_\RR(f)\circ \rho^{-1}.\label{5.14}
\end{eqnarray}
(\eqref{5.13} does not require assumption \eqref{5.1} or \eqref{5.2}).
As $G_\Z$ acts properly discontinuously on $M^{mar}_\mu(f_0)$,
$G^{smar}_\RR(f)$ and $G^{mar}_\RR(f)$ are finite.
(But this follows already from the finiteness of $R_f$ and (b).)
\end{theorem}

In the case of a quasihomogeneous singularity the group $R_f$ has a canonical lift
to $\RR^f$. It will be useful for the calculation of $R_f$.

\begin{theorem}\label{t5.9}
\cite[theorem 13.11]{He6}
Let $f\in \C[x_0,...,x_n]$ be a quasihomogeneous polynomial with an isolated singularity
at $0$ and weights $w_0,...,w_n\in\Q\cap(0,\frac{1}{2}]$ and weighted degree $1$.
Suppose that $w_0\leq ...\leq w_{n-1}<\frac{1}{2}$ (then $f\in \mmm^3$
if and only if $w_n<\frac{1}{2})$. Let $G_w$ be the algebraic group of 
quasihomogeneous coordinate changes, that means, those which respect $\C[x_0,...,x_n]$
and the grading by the weights $w_0,...,w_n$ on it. Then
\begin{eqnarray}\label{5.15}
R_f \cong \Stab_{G_w}(f).
\end{eqnarray}
\end{theorem}

Finally we need and we want to study period maps and 
Torelli type problems for singularities.

This story should start with the definition of the
Gau{\ss}-Manin connection and the Brieskorn lattice 
for an isolated hypersurface singularity. 
This had been developed in many papers of the second author,
and also much earlier by Brieskorn, K. Saito, G.-M. Greuel,
F. Pham, A. Varchenko, M. Saito and others. 

As we will build here on calculations done in \cite{He1}
and therefore never have to touch Brieskorn lattices explicitly,
we take here a formal point of view and refer to 
\cite{He1}, \cite{He2}, \cite{He4}, \cite{He6} and \cite{He7}
for the definitions of the following objects.

Any singularity $f$ comes equipped with a {\it Brieskorn lattice}
$H_0''(f)$. It is much richer than, but still comparable
to a Hodge structure of a closed K\"ahler manifold. 

After fixing a reference singularity $f_0$, 
a marked singularity $(f,\pm\rho)$  comes equipped with a 
{\it marked Brieskorn lattice} $BL(f,\pm\rho)$.
A classifying space $D_{BL}(f_0)$ for marked Brieskorn lattices was
constructed in \cite{He4}. It is especially a complex manifold,
and $G_\Z$ acts properly discontinuously on it.

\begin{theorem}\label{t5.10}
Fix one reference singularity $f_0$.

(a) There is a natural holomorphic period map
\begin{eqnarray}\label{5.16}
BL:M^{mar}_\mu(f_0)\to D_{BL}(f_0).
\end{eqnarray}
It is $G_\Z$-equivariant.

(b) \cite[theorem 12.8]{He6}
It is an immersion, here the reduced complex structure on 
$M^{mar}_\mu(f_0)$ is considered.
(The second author has also a proof that it is an immersion 
where the canonical complex structure on $M^{mar}_\mu(f_0)$ 
is considered, but the proof is not written).
\end{theorem}

The second main conjecture in \cite{He7} is part (a) of 
the following conjecture. 
Part (a) and part (b) are global Torelli type conjectures.

\begin{conjecture}\label{t5.11}
Fix one reference singularity $f_0$. 

(a) The period
map $BL:M^{mar}_\mu\to D_{BL}$ is injective.

(b) The period map
$LBL:M_\mu=M^{mar}_\mu/G_\Z\to D_{BL}/G_\Z$ is injective.

(c) For any singularity $f$ in the $\mu$-homotopy class
of $f_0$ and any marking $\rho$,
\begin{eqnarray}\label{5.17}
\textup{Stab}_{G_\Z}([(f,\pm\rho)])
=\textup{Stab}_{G_\Z}(BL([(f,\pm\rho)]))
\end{eqnarray}
(only $\subset$ and the finiteness of both groups are clear).
\end{conjecture}

The second author has a long-going project on Torelli type
conjectures. Already in \cite{He1}, part (b) was conjectured
and proved for all simple and unimodal singularities and 
almost all bimodal singularities 
(all except 3 subseries of the 8 bimodal series).
This was possible without the general construction of
$M_\mu$ and $D_{BL}$, which came later in \cite{He6}
and \cite{He4}. In the concrete cases considered in \cite{He1},
it is easy to identify a posteriori the spaces $M_\mu$
and $D_{BL}$. We will make use of that in the sections
\ref{c6} and \ref{c7}. Part (c) was conjectured in \cite{He6}.

The following lemma from \cite{He7} clarifies the logic
between the parts (a), (b) and (c) of conjecture \ref{t5.11}.

\begin{lemma}\label{t5.12}
In conjecture \ref{t5.11}, (a) $\iff$ (b) and (c).
\end{lemma}

Part (a) of conjecture \ref{t5.11} was proved in 
\cite{He7} for the simple and those 22 of the 28 exceptional
unimodal and bimodal singularities, where all eigenvalues
of the monodromy have multiplicity one.
Here it will be proved for the remaining unimodal and the
remaining exceptional bimodal singularities.

As part (b) of conjecture \ref{t5.11} was already proved
in all these cases in \cite{He1}, the main work in 
\cite[section 8]{He7} and here is the control of the
group $G_\Z$. This is carried out here in the 
sections \ref{c3} and \ref{c4}, and it is surprisingly difficult.

\section{$G^{mar}$, $M_\mu^{mar}$
and a strong Torelli result for the simple elliptic
and the hyperbolic singularities}\label{c6}

\noindent
The 1-parameter families of the hyperbolic singularities
of type $T_{pqr}$ ($p,q,r\in\N_{\geq 2}$, $p\geq q\geq r$,
$\kappa:=\frac{1}{p}+\frac{1}{q}+\frac{1}{r}<1$) have
as surface singularities the normal forms \cite{AGV1}
\begin{eqnarray}\label{6.1}
x^p+y^q+z^r+t\cdot xyz,\quad t\in X:=\C^*.
\end{eqnarray}
The 1-parameter families of the simple elliptic singularities
$T_{333}=\www E_6,T_{442}=\www E_7,T_{632}=\www E_8$
have as surface singularities different normal forms
\cite{SaK}.
The normal form $x^p+y^q+z^r+t\cdot xyz$ does in the case
of $T_{442}$ not contain representatives of all 
right equivalence classes, the class with $j$-invariant $j=1$ 
is missing \cite[1.11, Bem. (ii)]{SaK}. 
Therefore we work in the following
also with the Legendre normal forms
\begin{eqnarray}
T_{333}:&& y(y-x)(y-\lambda x)-xz^2,
\quad  t\in X:=\C-\{0,1\},\nonumber\\
T_{442}:&& yx(y-x)(y-\lambda x)+z^2,
\quad  t\in X:=\C-\{0,1\},\label{6.2}\\
T_{632}:&& y(y-x^2)(y-\lambda x^2)+z^2,
\quad  t\in X:=\C-\{0,1\}.\nonumber
\end{eqnarray}
They contain representatives of all right equivalence classes.
Let $X^{univ}$ denote the universal covering of $X$,
so $X^{univ}=\C$ if $\kappa<1$ and $X^{univ}=\H$ if $\kappa=1$.

\begin{theorem}\label{t6.1}
(a) For the simple elliptic singularities and the hyperbolic 
singularities in any number of variables, the space $M^{mar}_\mu$
of right equivalence classes of marked singularities is
$M^{mar}_\mu\cong X^{univ}$, so it is connected, 
and thus $G^{mar}=G_\Z$. The period map 
$BL:M^{mar}_\mu\to D_{BL}$ is an isomorphism, so the 
strong global Torelli conjecture \ref{t5.11} (a) is true.

(b) Now consider the singularities of type $T_{pqr}$
as curve singularities if $r=2$ and as surface singularities
if $r\geq 3$. Then
\begin{eqnarray}\label{6.3}
G_\Z=G^{mar}=G^{smar}\times\{\pm\id\},\quad\textup{equivalently: }
-\id\notin G^{smar}.
\end{eqnarray}
The subgroup of $G^{smar}$, which acts trivially on
$M^{mar}_\mu$, is the kernel of the surjective map
\begin{eqnarray}\label{6.4}
G^{smar}\to \Aut(Ml(f_0)_{1,\Z},L)/\{\pm\id\}.
\end{eqnarray}
It is equal to $\rho\circ (\RR^f)_{hom}\circ\rho^{-1}$
for a generic $[(f,\rho)]\in M^{mar}_\mu$. 
Its size is 54, 16 and 6 for $T_{333}$, $T_{442}$ and $T_{632}$.
\end{theorem}

{\bf Proof:}
(a) The proof uses two Torelli type results from \cite{He1}.

We choose a marked reference singularity $[(f_0,\pm\id)]$
in $X^{univ}$, then all elements of $X^{univ}$ 
become marked singularities, because $X^{univ}$ is 
simply connected.
Then the period map $X^{univ}\to D_{BL}$ is well defined. 
The first Torelli type result from \cite{He1} is that
this map is an isomorphism. 

Therefore the marked Brieskorn lattices of the marked singularities
in $X^{univ}$ are all different.
Therefore the marked singularities in $X^{univ}$ are
all not right equivalent. This gives an embedding
$X^{univ}\hookrightarrow M^{mar}_\mu(f_0)^0$. 

On the other hand, we have the period map 
$BL: M^{mar}_\mu(f_0)^0\to D_{BL}$, which is an immersion.
As it restricts to the isomorphism $X^{univ}\to D_{BL}$, finally
$X^{univ}=M^{mar}_\mu(f_0)^0$.

For part (a) it rests to show that $G^{mar}=G_\Z$.
Then $M^{mar}_\mu$ is connected and $M^{mar}_\mu=X^{univ}$,
and $BL:M^{mar}_\mu\to D_{BL}$ is an isomorphism.

We have to look closer at $D_{BL}$ and the action of $G_\Z$
on it. In the case $\kappa<1$, 
\begin{eqnarray}\label{6.5}
D_{BL}\cong \{V\subset Ml(f_0)_1\, |\, \dim V=1, 
V\neq \ker(M_h-\id)\}.
\end{eqnarray}
In the case $\kappa =1$,
\begin{eqnarray}\label{6.6}
D_{BL}\cong \textup{one component of }
\{V\subset Ml(f_0)_1\, |\, \dim V=1, 
V\neq \oooo{V}\}.
\end{eqnarray}
In both cases the group $\Aut(Ml(f_0)_{1,\Z},L)/\{\pm \id\}$
acts faithfully on $D_{BL}$.

In both cases the second Torelli type result from \cite{He1}
which we need is that the period map
\begin{eqnarray}\label{6.7}
X^{univ}/\sim_R\ \to D_{BL}/ \Aut(Ml(f_0)_{1,\Z},L)
\end{eqnarray}
is an isomorphism. Here $\sim_R$ denotes right equivalence
for unmarked singularities.

Because of $X^{univ}=M^{mar}_\mu(f_0)^0$, 
\begin{eqnarray}\label{6.8}
X^{univ}/\sim_R = M^{mar}_\mu(f_0)^0/G^{mar}.
\end{eqnarray}
The isomorphism \eqref{6.7} and the isomorphism
$M^{mar}_\mu(f_0)^0=X^{univ}\to D_{BL}$ show that the map
\begin{eqnarray}\label{6.9}
G^{mar}\to \Aut(Ml(f_0)_{1,\Z},L)
\end{eqnarray}
is surjective. This completes the proof of part (b) of theorem \ref{t3.1}.

It also shows that for proving $G^{mar}=G_\Z$, it is sufficient
to show that the kernels of the maps to
$\Aut(Ml(f_0)_{1,\Z},L)/\{\pm\id\}$ coincide.
The kernel of the map $G_\Z \to \Aut(Ml(f_0)_{1,\Z},L)/\{\pm\id\}$
was determined in theorem \ref{t3.1}. In fact, this is 
the only part of theorem \ref{t3.1} which we need here. 
It consists of those elements of $(U_1\rtimes U_2)\times\{\pm\id\}$ 
in \eqref{3.6} for which $\delta=0$, so it is isomorphic
to the group
\begin{eqnarray}
&&\hspace*{-0.7cm}\left(\{(\alpha,\beta,\gamma)\in\Z_p\times\Z_q\times\Z_r\, |\, 
\frac{\alpha}{p}+\frac{\beta}{q}+\frac{\gamma}{r}\equiv 0
\mod 1\}\rtimes U_2\right)\times\{\pm\id\}\nonumber\\
&&=: (U_1^0\rtimes U_2)\times\{\pm\id\}.\label{6.10}
\end{eqnarray}

The kernel of the map $G^{mar}\to\Aut(Ml(f_0)_{1,\Z},L)/\{\pm\id\}$
is the subgroup of $G^{mar}$ which acts trivially on 
$M^{mar}_\mu(f_0)^0$. It is the isotropy group in $G^{mar}$
of a generic point $[(f,\pm\rho)]\in M^{mar}_\mu(f_0)^0$.
So by theorem \ref{t5.8} (g) it is the group
\begin{eqnarray}\label{6.11}
\rho\circ G^{mar}_\RR(f)\circ \rho^{-1}
=\rho\circ \{\pm\varphi_{hom}\, |\, \varphi\in\RR^f\}\circ\rho^{-1}.
\end{eqnarray}

As $f$ is generic, we can and will use now the normal form
$f=x^p+y^q+z^r+t\cdot xyz$, also in the case $\kappa=1$.
The following coordinate changes generate a finite subgroup
$S\subset \RR^f$ 
\begin{eqnarray}\label{6.12}
\varphi^{\alpha,\beta,\gamma}: (x,y,z)&\mapsto&
(e^{2\pi i\alpha/p}x,e^{2\pi i\beta/q}y,e^{2\pi i\gamma/r}z)\\
&&ß\textup{with }\frac{\alpha}{p}+\frac{\beta}{q}+\frac{\gamma}{r}
\equiv 0\mod 1,\nonumber\\
\varphi^{1,2}: (x,y,z)&\mapsto& 
(y,x,z)\quad\textup{if }p=q,\nonumber\\
\varphi^{2,3}: (x,y,z)&\mapsto& 
(x,z,y)\quad\textup{if }q=r,\nonumber\\
\varphi^{minus}: (x,y,z)&\mapsto& 
(x,y,-z-txy)\quad\textup{if }r=2.\nonumber
\end{eqnarray}
$\varphi^{minus}$ has order 2 and commutes with the other
coordinate changes ($q=r=2$ is impossible because of $\kappa\leq 1$).
The group $S$ is isomorphic (as an abstract group) 
to $U_1^0\rtimes U_2$ if $r\geq 3$
and to $(U_1^0\rtimes U_2)\times\{\pm\id\}$ if $r=2$.
The map to 1-jets of coordinate changes is injective,
\begin{eqnarray}\label{6.13}
S\stackrel{\cong}{\longrightarrow} 
j_1S\subset j_1\RR^f\subset j_1\RR.
\end{eqnarray}
Now we have to treat the cases $r\geq 3$ and $r=2$ separately.

{\bf The case $r\geq 3$:}
Then $j_1\RR^f$ is finite and isomorphic to $R_f$, the map 
$$()_{hom}:R_f\to G_\Z(f)=\rho^{-1}\circ G_\Z\circ \rho$$
is injective, and the image $G^{smar}_\RR(f)$ 
does not contain $-\id$ (theorem \ref{t5.8}).
Therefore then $S\cong (S)_{hom}\subset G_\Z(f)$ and 
$-\id\notin (S)_{hom}$. Thus the group $(S)_{hom}\times\{\pm\id\}$
is isomorphic to $(U_1^0\rtimes U_2)\times\{\pm\id\}$.
Now it is clear that the group in \eqref{6.11} is at least
as big as the group in \eqref{6.10}. But it cannot be bigger.
So they are of equal size. This implies $G^{mar}=G_\Z$.

{\bf The case $r=2$:}
We claim that the map $S\to (S)_{hom}$ is injective.
If this is true then 
$(S)_{hom}\cong (U_1^0\rtimes U_2)\times\{\pm\id\}$,
and this is of equal size as the group in \eqref{6.10}.
Then again the group in \eqref{6.11} is at least as big as
the group in \eqref{6.10}, but it cannot be bigger.
So they are of equal size. This implies $G^{mar}=G_\Z$.

It rests to prove the claim. For this we consider the
curve singularity 
\begin{eqnarray}\label{6.14}
g:=x^p+y^q-\frac{1}{4}tx^2y^2.
\end{eqnarray}
Then
\begin{eqnarray}
g+z^2 &=&f\circ\psi\quad\textup{with }
\psi(x,y,z)=(x,y,z-\frac{1}{2}txy),\nonumber\\
\RR^{g+z^2} &=& \psi^{-1}\circ\RR^f\circ\psi,\label{6.15}\\
\psi^{-1}\circ \varphi^{\alpha,\beta,\gamma}\circ \psi
&=& \varphi^{\alpha,\beta,\gamma},\nonumber\\
\psi^{-1}\circ \varphi^{1,2}\circ \psi
&=& \varphi^{1,2},\quad\textup{if }p=q,\nonumber\\
\psi^{-1}\circ \varphi^{minus}\circ \psi
&=& ((x,y,z)\mapsto (x,y,-z)).\nonumber
\end{eqnarray}
($q=r=2$ is impossible because of $\kappa\leq 1$). The subgroup 
\begin{eqnarray}\label{6.16}
S^{curve}:= \{\varphi^{\alpha,\beta,\gamma}\circ 
(\varphi^{minus})^{-\gamma}\, |\, (\alpha,\beta,\gamma)\in U_1^0\}
\rtimes U_2
\end{eqnarray}
has index 2 in $S$, its conjugate $\psi^{-1}\circ S^{curve}\circ\psi$
restricts to $\RR^g$, and it maps injectively to
$j_1\RR^g\cong R_g$. By theorem \ref{t5.8} (c) the map
$S^{curve}\to (S^{curve})_{hom}$ is injective, and $-\id$
is not in the image. But $(\varphi^{minus})_{hom}=-\id$.
This proves the claim.

\medskip
(b) By theorem \ref{t5.5} (c), the projection
$M_\mu^{smar}\to M^{mar}_\mu$ is a twofold covering, and 
$-\id$ exchanges the two sheets of the covering.
Because of $M^{mar}_\mu=\C$ if $\kappa<1$ and $M^{mar}_\mu=\H$
if $\kappa=1$, $M^{smar}_\mu$ has two components.
Therefore $-\id\notin G^{smar}$ and 
$G_\Z=G^{mar}=G^{smar}\times\{\pm\id\}$.

The statements right before and after \eqref{6.4} were already
proved and used in the proof of part (a). 

The group $\rho\circ (\RR^f)_{hom}\circ \rho^{-1}$ for a generic 
$[(f,\rho)]\in M^{mar}_\mu$ has size 54, 16 and 6 for
$T_{333}, T_{442}$ and $T_{632}$, because it is isomorphic
to an index 2 subgroup of the group in \eqref{6.10},
and that group has 108, 32 and 12 elements in the cases
$T_{333}, T_{442}$ and $T_{632}$ \hfill $\Box$

\section{$G^{mar}$, $M_\mu^{mar}$
and a strong Torelli result for 6 of the 28 exceptional
unimodal and bimodal singularities}\label{c7}

\noindent
Normal forms for the 1-parameter families of the exceptional
unimodal and bimodal singularities of types 
$Z_{12}, Q_{12}, U_{12}, Z_{18}, Q_{16}$ and $U_{16}$
in the minimal number of variables are as follows
\cite{AGV1}. Here $Z_{12}$ and $Z_{18}$ are curve singularities,
$Q_{12}, U_{12}, Q_{16}$ and $U_{16}$ are surface singularities.
The singularity for the parameter $t=0$ is
quasihomogeneous, the others are semiquasihomogeneous.
The space of the parameter $t=t_1$ or $t=(t_1,t_2)$ is 
$X=\C^{\textup{mod}(f_0)}$.
The weights $w=(w_x,w_y)$ respectively 
$w=(w_x,w_y,w_z)$ are normalized such that $\deg_wf_0=1$.

\begin{eqnarray}\label{7.1}
\begin{array}{cllcl}
 & & \textup{normal form} & \textup{mod}(f_0) &
\textup{weights }\\ \hline
Z_{12}:& f_t &=x^3y+xy^4+tx^2y^3& 1 
& (\frac{3}{11},\frac{2}{11})\\
Q_{12}:& f_t &=x^3+y^5+yz^2+txy^4& 1
& (\frac{1}{3},\frac{1}{5},\frac{2}{5})\\
U_{12}:& f_t &=x^3+y^3+z^4+txyz^2&1
& (\frac{1}{3},\frac{1}{3},\frac{1}{4})\\
Z_{18}:& f_t &=x^3y+xy^6 + (t_1+t_2y)y^9& 2
& (\frac{5}{17},\frac{2}{17})\\
Q_{16}:& f_t &=x^3+y^7+yz^2 + (t_1+t_2y)xy^5& 2
& (\frac{1}{3},\frac{1}{7},\frac{3}{7})\\
U_{16}:& f_t &=x^3+xz^2+y^5 + (t_1+t_2y)x^2y^2& 2
& (\frac{1}{3},\frac{1}{5},\frac{1}{3})
\end{array}
\end{eqnarray}

The normal form of the quasihomogeneous singularity of type
$Q_{12}$, $Q_{16}$, $U_{12}$ and $U_{16}$ is a sum
of an $A_l$-singularity in one variable and a 
$D_{2m}$-singularity in two variables with $(l,2m)$
as in table \eqref{7.2}=\eqref{4.20}.

\begin{eqnarray}\label{7.2}
\begin{array}{l|l|l|l|l}
 &Q_{12}&Q_{16}&U_{12}&U_{16}\\ \hline
(l,2m)&(2,6)&(2,8)&(3,4)&(4,4)\\
 & A_2\otimes D_6 & A_2\otimes D_8 &A_3\otimes D_4 & A_4\otimes D_4
\end{array}
\end{eqnarray}

The rest of this section is devoted to the proof of the
following theorem.

\begin{theorem}\label{t7.1}
(a) For the 6 families of exceptional unimodal and bimodal
singularities of types $Z_{12}, Q_{12}, U_{12}, Z_{18}, Q_{16}$
and $U_{16}$ in any number of variables, the space $M^{mar}_\mu$
of right equivalence classes of marked singularities is
$M^{mar}_\mu\cong X=\C^{\textup{mod}(f_0)}$, so it is connected, 
and thus $G^{mar}=G_\Z$. The period map 
$BL:M^{mar}_\mu\to D_{BL}$ is an isomorphism, so the 
strong global Torelli conjecture \ref{t5.11} (a) is true.

(b) Now consider the singularities of type $Z_{12}$ and $Z_{18}$
as curve singularities and the singularities of types
$Q_{12}, U_{12}, Q_{16}$ and $U_{16}$ as surface singularities.
Then their multiplicities are $\geq 3$. Then
\begin{eqnarray}\label{7.3}
G_\Z=G^{mar}=G^{smar}\times\{\pm\id\},\quad\textup{equivalently: }
-\id\notin G^{smar}.
\end{eqnarray}
\end{theorem}

{\bf Proof:}
(a) The proof is similar to the proof of theorem \ref{t6.1},
but simpler. We need only the first of the two Torelli type
results from \cite{He1}, which were used in the proof of 
theorem \ref{t6.1}. 

We choose as marked reference singularity the 
quasihomogeneous singularity with trivial marking $[(f_0,\pm\id)]$
in $X$. Then all elements of $X$ 
become marked singularities, because $X$ is 
simply connected.
Then the period map $X\to D_{BL}$ is well defined. 
A Torelli type result from \cite{He1} says that
this map is an isomorphism. It is in fact easy, 
because the singularities here are semiquasihomogeneous 
and only $f_0$ is quasihomogeneous. That makes the calculations
easy.

Therefore the marked Brieskorn lattices of the marked singularities
$X$ are all different.
Therefore the marked singularities in $X$ are
all not right equivalent. This gives an embedding
$X\hookrightarrow M^{mar}_\mu(f_0)^0$. 

On the other hand, we have the period map 
$BL: M^{mar}_\mu(f_0)^0\to D_{BL}$, which is an immersion.
As it restricts to the isomorphism $X\to D_{BL}$, finally
$X=M^{mar}_\mu(f_0)^0$. 

For part (a) it rests to show that $G^{mar}=G_\Z$.
Then $M^{mar}_\mu$ is connected and $M^{mar}_\mu=X$,
and $BL:M^{mar}_\mu\to D_{BL}$ is an isomorphism.

The weights of the deformation parameter(s) $t_1$ (and $t_2$)
equip the parameter space $X=M^{mar}_\mu$ with a good 
$\C^*$-action. It commutes with the action of $G_\Z$. 
This gives the first equality in 
\begin{eqnarray}\label{7.4}
G^{mar}=\textup{Stab}_{G_\Z}([(f_0,\pm\id)])=G^{mar}_\RR(f_0).
\end{eqnarray}
The second equality is part of theorem \ref{t5.8} (g).

On the other hand, by theorem \ref{t5.8} (d),
$G^{smar}_\RR(f_0)\times\{\pm\id\}=G^{mar}_\RR(f_0)$.
But because $f_0$ is a quasihomogeneous singularity of degree
$\geq 3$, the group $G^{smar}_\RR(f_0)$ can be calculated easily via 
$\Stab_{G_w}(f_0)$, see theorem \ref{t5.8} (c) and 
theorem \ref{t5.9}: 
$\textup{Stab}_{G_w}(f_0)\stackrel{\cong}{\longrightarrow}
G^{smar}_\RR(f_0)$.
Therefore it is sufficient to show that
$\textup{Stab}_{G_w}(f_0)$ 
has half as many elements as the group $G_\Z$. 
We postpone its proof. If it holds, then
\begin{eqnarray}\label{7.5}
G_\Z=G^{mar}=G^{mar}_\RR(f_0)=G^{smar}_\RR(f_0)\times\{\pm\id\}
\end{eqnarray}
follows, and part (a) of the theorem is proved.
For part (b), the same argument as in the proof of theorem
\ref{t6.1} (b) works: $M^{smar}_\mu$ is a twofold covering
of $M^{mar}_\mu$, and the two sheets are exchanged by the
action of $-\id$. As $X=\C^{\textup{mod}(f_0)}$, $M^{smar}_\mu$
has two components, and $-\id\notin M^{smar}_\mu$.

In theorem \ref{t4.1} it was shown that
$G_\Z$ is $G_\Z=\{\pm M_h^k\, |\, k\in\Z\}\times U$
with $U$ as in table \eqref{4.1}=\eqref{7.6}.
\begin{eqnarray}\label{7.6}
\begin{array}{l|l|l|l|l|l|l}
 & Z_{12} & Q_{12} & U_{12} & Z_{18} & Q_{16} & U_{16}\\ \hline
U\cong & \{\id\} & S_2 & S_3 & \{\id\} & S_2 & S_3
\end{array}
\end{eqnarray}
Now we compare $\textup{Stab}_{G_w}(f_0)$.
It is sufficient to find enough elements so that the resulting
group has half as many elements as $G_\Z$.

{\bf The cases $Z_{12}$ and $Z_{18}$:}
Then 
\begin{eqnarray}\label{7.7} 
\varphi_1:(x,y)\mapsto (e^{2\pi i w_x}x,e^{2\pi iw_y}y)
\quad\textup{satisfies}\quad (\varphi_1)_{hom}=M_h.
\end{eqnarray}
This is already sufficient. Here $G^{smar}=G^{smar}_\RR
=\{M_h^k\, |\, k\in\Z\}.$

{\bf The cases $Q_{12},Q_{16},U_{12},U_{16}$:}
Here it is convenient to make use of the decomposition
of the singularity $f_0$ into a sum of an $A_l$ singularity $g_0$ 
in one variable and a $D_{2m}$ singularity $h_0$ in two variables.
In all four cases the weight system $w'$ of the $A_l$ singularity
and the weight system $w''$ of the $D_{2m}$ singularity 
have denominators $l+1$ and $2m-1$ with $\gcd(l+1,2m-1)=1$.
Therefore 
\begin{eqnarray}\label{7.8}
\textup{Stab}_{G_w}(f_0)
&=&\textup{Stab}_{G_{w'}}(g_0)\times \textup{Stab}_{G_{w''}}(h_0)\\
&\cong& \Z_{l+1}\times \left\{
\begin{array}{ll} \Z_{2m-1}\times S_2&\textup{ if }m\geq 3\\
\Z_3\times S_3&\textup{ if }m=2.\end{array}\right.  \nonumber
\end{eqnarray}
In all four cases this group has half as many elements as $G_\Z$.
\hfill$\Box$

\section{More on $G_\Z$ for the simple elliptic singularities}\label{c8}

\noindent
This section is motivated by the paper \cite{MS} of Milanov and Shen.
They consider the 1-parameter families 
\begin{eqnarray}\label{8.1}
x^p+y^q+z^r+t\cdot xyz, \quad t\in\Sigma\subset\C
\end{eqnarray}
of the simple elliptic singularities $T_{pqr}$ with
$(p,q,r)\in\{(3,3,3),(4,4,2),(6,3,2)\}$ which had also been
used above in section \ref{c6}. Here $\Sigma\subset\C$
is the complement of the finite set of parameters where the function
in \eqref{8.1} has a non-isolated singularity. 
Remark that now $\chi:=\lcm(p,q,r)=p$.

In \cite{MS} the groups of the transversal monodromies
of these three families are studied,
more precisely, the natural representations
\begin{eqnarray}\label{8.2}
\rho:\pi_1(\Sigma)&\to& G_\Z,\\
\rho_1:\pi_1(\Sigma)&\to&\Aut(Ml(f)_{1,\Z},L),\nonumber\\
\rho_{\neq 1}:\pi_1(\Sigma)&\to&\Aut(Ml(f)_{\neq 1,\Z},L),\nonumber\\
\oooo{\rho}_{\neq 1}:\pi_1(\Sigma)&\to&\Aut(Ml(f)_{\neq 1,\Z},L)/\langle M_h\rangle.
\nonumber
\end{eqnarray}
By explicit computations they show
\begin{eqnarray}\label{8.3}
\ker(\rho_1)\subset\ker(\oooo{\rho}_{\neq 1}).
\end{eqnarray}
They ask about a conceptual explanation of \eqref{8.3} and whether this might
be true for other normal forms, i.e. other natural 1-parameter families
of the simple elliptic singularities.
Because of \eqref{8.3}, there is an induced representation
\begin{eqnarray}\label{8.4}
\rho_W:\Imm(\rho_1)\to \Aut(Ml(f)_{\neq 1},L)/\langle M_h\rangle.
\end{eqnarray}
Then $\ker(\rho_W)\subset \Imm(\rho_1)\subset\Aut(Ml(f)_{1,\Z},L)$. 
One main result of \cite{MS} is the following.

\begin{theorem}\label{t8.1}
In all three cases, 
under an isomorphism $\Aut(Ml(f)_{1,\Z},L)\cong SL(2,\Z)$ as in \eqref{3.2}, 
the subgroup $\ker(\rho_W)$ is isomorphic to the 
principal congruence subgroup $\Gamma(p)$.
\end{theorem}

Here $\Gamma(N):=\{A\in SL(2,\Z)\, |\, A\equiv {\bf 1}_2\mod N\}$
(not $A\equiv\pm{\bf 1}_2\mod N$) is the principal congruence subgroup
of level $N$ of $SL(2,\Z)$.

In the following, we give
results which complement \eqref{8.3} and theorem \ref{t8.1}.
We consider not a special 1-parameter family, but the biggest possible
family of quasihomogeneous simple elliptic singularities.
Define $(w_x,w_y,w_z):=(\frac{1}{p},\frac{1}{q},\frac{1}{r})$, define
for any monomial its weighted degree $\deg _w(x^\alpha y^\beta z^\gamma)
:=\alpha w_x+\beta w_y+\gamma w_z$, and define
\begin{eqnarray}\label{8.5}
&&\C[x,y,z]_1:=\langle x^\alpha y^\beta z^\gamma\, |\, 
\deg_w(x^\alpha y^\beta z^\gamma)=1\rangle_\C,\\
&&R:=\{f\in\C[x,y,z]_1\, |\, f\textup{ has an isolated singularity at 0}\}.
\nonumber
\end{eqnarray}
$R$ is the complement of a hypersurface in the vector space $\C[x,y,z]_1$.
The transversal monodromies of the family of singularities parametrized
by $R$ give the natural representation $\sigma$, the other representations
are induced,
\begin{eqnarray}\label{8.6}
\sigma:\pi_1(R)&\to& G_\Z,\\
\sigma_1:\pi_1(R)&\to&\Aut(Ml(f)_{1,\Z},L),\nonumber\\
\sigma_{\neq 1}:\pi_1(R)&\to&\Aut(Ml(f)_{\neq 1,\Z},L),\nonumber\\
\oooo{\sigma}_{\neq 1}:\pi_1(R)&\to&
\Aut(Ml(f)_{\neq 1,\Z},L)/\langle M_h\rangle.\nonumber
\end{eqnarray}
By the definition of $G^{smar}$ in \cite[definition 3.1]{He7},
$G^{smar}=\Imm(\sigma)$. Theorem \ref{t6.1} tells that the monodromy group
$\Imm(\sigma)$ is as large as possible (up to $\pm\id$ in the case $T_{333}$).

\begin{theorem}\label{t8.2}
$\Imm(\sigma)=G_\Z$ in the cases $T_{442}$ and $T_{632}$,
and $G_\Z=\Imm(\sigma)\times\{\pm\id\}$ in the case $T_{333}$
(here it is important that the surface singularities are considered).
\end{theorem}

The explicit information on $G_\Z$ in theorem \ref{t3.1} allows the following
conclusion.

\begin{corollary}\label{t8.3}
The analogue $\ker(\sigma_1)\subset\ker(\oooo{\sigma}_{\neq 1})$ 
of \eqref{8.3}
does not hold in the cases $T_{333}$ and $T_{442}$. It holds in the
case $T_{632}$, and there the analogue of \eqref{8.3} holds
for any $\mu$-constant family.
\end{corollary}

{\bf Proof:} By theorem \ref{t3.1} (c), 
\begin{eqnarray*}
&&\{g\in G_\Z\, |\, g|_{Ml(f)_1}=\id\}\\
&=& \{\id|_{Ml(f)_1}\times (M_h|_{Ml^{(1)}_\C})^\alpha
\times (M_h|_{Ml^{(2)}_\C})^\beta
\times (M_h|_{Ml^{(3)}_\C})^\gamma\, | \\
&&(\alpha,\beta,\gamma)\in\Z_p\times\Z_q\times \Z_r\textup{ with }
\frac{\alpha}{p}+\frac{\beta}{q}+\frac{\gamma}{r}\equiv0\mod 1\}\times U_2.
\end{eqnarray*}
In the cases $T_{442}$ and $T_{632}$ this is isomorphic to 
$\ker(\sigma_1)/\ker(\sigma)$, 
in the case $T_{333}$ $\ker(\sigma_1)/\ker(\sigma)$ 
is isomorphic to this group or a subgroup of this group
of index 2. In the cases $T_{442}$ and $T_{333}$, already the
factor $U_2$ in this group is an obstruction to the analogue of \eqref{8.3}.

In the case $T_{632}$, $U_2$ is trivial and
\begin{eqnarray}\label{8.8}
&&\{g\in G_\Z\, |\, g|_{Ml(f)_1}=\id\}=\{M_h^\alpha\, |\, \alpha\in\Z\}
\end{eqnarray}
(because above $\alpha$ determines $\beta$ and $\gamma$ uniquely
in the case $(p,q,r)=(6,3,2)$).
Therefore the analogue of \eqref{8.3} holds in the case $T_{632}$
for the family parametrized by $R$ and for any subfamily. \hfill$\Box$

\bigskip
\eqref{8.3} is used in \cite{MS} in order to define $\rho_W$ and 
the group $\ker(\rho_W)$.
But because $M_h|_{Ml(f)_1}=\id$, it is obvious that the group
$\ker(\rho_W)$ coincides with 
\begin{eqnarray}\label{8.9}
\{g|_{Ml(f)_{1,\Z}}\, |\ g\in \Imm(\sigma),g|_{Ml(f)_{\neq 1,\Z}}=\id\}
\subset\Aut(Ml(f)_{1,\Z},L).
\end{eqnarray} 
And the analogue of this group 
can be defined for any $\mu$-constant family, whether or not it satisfies
the analogue of \eqref{8.3}. Our main result in this section is 
theorem \ref{t8.4}. Our proof uses theorem \ref{t3.1}. 
A different proof of theorem \ref{t8.4} was given by Kluitmann in
\cite[III 2.4 Satz, page 66]{Kl}. 
Theorem \ref{t8.4} shows that the group $\Gamma(p)$ turns up naturally
within the maximal possible $\mu$-constant family, which is 
parametrized by $R$,
and it shows the part $\ker(\rho_W)\subset \Gamma(p)$ of 
the equality $\ker(\rho_W)= \Gamma(p)$ in theorem \ref{t8.1}.

\begin{theorem}\label{t8.4}
In all three cases, under an isomorphism 
$\Aut(Ml(f)_{1,\Z},L)\cong SL(2,\Z)$ as in \eqref{3.2}, 
the subgroup 
$$\{g|_{Ml(f)_{1,\Z}}\, |\ g\in G_\Z,g|_{Ml(f)_{\neq 1,\Z}}=\id\}
\subset\Aut(Ml(f)_{1,\Z},L).$$ 
is isomorphic to the principal congruence subgroup $\Gamma(p)$.
\end{theorem}

{\bf Proof:} We use the notations and objects in the proof of theorem
\ref{t3.1}. $Ml^{(1)}_\C$ was defined in \eqref{3.15}. Define
\begin{eqnarray}\label{8.10}
Ml^{(1)}_\Z:=Ml^{(1)}_\C\cap Ml(f),
\end{eqnarray}
and analogously $Ml^{(2)}_\Z$ and $Ml^{(3)}_\Z$. Then 
$$Ml^{(1)}_\Z=(M_h-\id)(Ml^{[1]}_\Z).$$
Because of \eqref{3.16}, this image is generated as a $\Z$-lattice by
\begin{eqnarray}\label{8.11}
\delta_2-(\delta_1+\www b_1), \delta_3-\delta_2,...,\delta_{p-1}-\delta_{p-2},
-(\delta_1+...+\delta_{p-1})-\delta_{p-1},
\end{eqnarray}
respectively by
\begin{eqnarray}\label{8.12}
\delta_1+(p-1)\delta_2,\delta_3-\delta_2,...,\delta_{p-1}-\delta_{p-2},
p\delta_2-\delta_{\mu-1}+\delta_\mu.
\end{eqnarray}
$Ml^{(2)}_\Z$ and $Ml^{(3)}_\Z$ (if $r\geq 3$) are generated by the 
analogous elements. If $r=2$ then 
\begin{eqnarray}\label{8.13}
Ml^{(3)}_\Z =\Z\cdot (2\delta_{\mu-2}-\delta_{\mu-1}+\delta_\mu).
\end{eqnarray}
In any case, the sum $Ml^{(1)}_\Z\oplus Ml^{(2)}_\Z\oplus Ml^{(3)}_\Z$
is a sublattice of finite index of the primitive sublattice 
$Ml(f)_{\neq 1,\Z}$ in $Ml(f)$. Observe that $q|p$ and $r|p$ in all three
cases. The lattice $Ml(f)_{\neq 1,\Z}$ is generated by
\begin{eqnarray}\label{8.14}
\delta_1+(p-1)\delta_2,\delta_3-\delta_2,...,\delta_{p-1}-\delta_{p-2},
p\delta_2-\delta_{\mu-1}+\delta_\mu,\\
\delta_p+(q-1)\delta_{p+1},\delta_{p+2}-\delta_{p+1},...,\delta_{p+q-2}-\delta_{p+q-3},
\frac{p}{q}\delta_2-\delta_{p+1},\nonumber \\
\delta_{p+q-1}+(r-1)\delta_{p+q},\delta_{p+q+1}-\delta_{p+q},...,
\delta_{\mu-2}-\delta_{\mu-3},\frac{p}{r}\delta_2-\delta_{p+q},\nonumber
\end{eqnarray}
if $r\geq 3$. If $r=2$ then the third line has to be replaced by
$$\frac{p}{r}\delta_2-\delta_{\mu-2}.$$
In any case, one sees 
\begin{eqnarray}\label{8.15}
Ml(f)=Ml(f)_{\neq 1,\Z}\oplus\Z\cdot\delta_2\oplus\Z\cdot\delta_\mu.
\end{eqnarray}
One also calculates 
\begin{eqnarray}\label{8.16}
\www b_1&=&\delta_{\mu-1}-\delta_\mu=\gamma_1+p\delta_2\\
\textup{with}\quad\gamma_1&:=&-(p\delta_2-\delta_{\mu-1}+\delta_\mu)
\in Ml(f)_{\neq 1,\Z},\nonumber\\
\www b_2 &=&\gamma_2+p\delta_\mu,\label{8.17}\\
\textup{with}\quad\gamma_2&:=&
\sum_{i=1}^{p-1}(p-i)\delta_i
+\sum_{i=1}^{q-1}\frac{p}{q}(q-i)\delta_{p-1+i}\nonumber\\
&&+\sum_{i=1}^{r-1}\frac{p}{r}(r-i)\delta_{p+q-2+i}
+p\delta_{\mu-1}-p\delta_\mu\in Ml(f)_{\neq 1,\Z}.\nonumber
\end{eqnarray}
One sees also
\begin{eqnarray}\label{8.18}
Ml(f)\cap(\Q\cdot \gamma_1+\Q\cdot\gamma_2)
=\Z\cdot\gamma_1\oplus \Z\cdot\gamma_2.
\end{eqnarray}
For any matrix $\begin{pmatrix}a&b\\c&d\end{pmatrix}\in SL(2,\Z)$
define the automorphism $f:Ml(f)_\Q\to Ml(f)_\Q$ by
\begin{eqnarray}\label{8.19}
f(\www b_1):=a\www b_1+c\www b_2,\quad f(\www b_2):=b\www b_1+d\www b_2,\quad
f|_{Ml(f)_{\neq 1,\Z}}:=\id.
\end{eqnarray}
It respects $L$ because the decomposition $Ml(f)_\Q=Ml(f)_{1,\Q}\oplus
Ml(f)_{\neq 1,\Q}$ is left and right orthogonal with respect to $L$.
It is an automorphism of $Ml(f)$ if and only if $f(\delta_2)$ and $f(\delta_\mu)$
are in $Ml(f)$. One calculates
\begin{eqnarray*}
f(\delta_2)&=&a\delta_2+c\delta_\mu+\frac{1}{p}(a-1)\gamma_1+\frac{1}{p}c\gamma_2,\\
f(\delta_\mu)&=&b\delta_2+d\delta_\mu+\frac{1}{p}b\gamma_1+\frac{1}{p}(d-1)\gamma_2.
\end{eqnarray*}
In view of \eqref{8.18} this shows
\begin{eqnarray*}
f\in G_\Z \iff \begin{pmatrix}a&b\\c&d\end{pmatrix}\equiv {\bf 1}_2\mod p
\stackrel{\textup{def.}}{\iff} \begin{pmatrix}a&b\\c&d\end{pmatrix}\in \Gamma(p).\hspace*{1cm}\Box
\end{eqnarray*}

\end{document}